\documentclass[10pt]{article}
\usepackage{amssymb, amsmath}
\usepackage{float,epsfig}
\usepackage{rotating}

\begin{document}

\newtheorem{theorem}{\bf Theorem}[section]
\newtheorem{proposition}[theorem]{\bf Proposition}
\newtheorem{definition}[theorem]{\bf Definition}
\newtheorem{corollary}[theorem]{\bf Corollary}
\newtheorem{example}[theorem]{\bf Example}
\newtheorem{exam}[theorem]{\bf Example}
\newtheorem{remark}[theorem]{\bf Remark}
\newtheorem{lemma}[theorem]{\bf Lemma}
\newcommand{\nrm}[1]{|\!|\!| {#1} |\!|\!|}

\newcommand{\ba}{\begin{array}}
\newcommand{\ea}{\end{array}}
\newcommand{\von}{\vskip 1ex}
\newcommand{\vone}{\vskip 2ex}
\newcommand{\vtwo}{\vskip 4ex}
\newcommand{\dm}[1]{ {\displaystyle{#1} } }

\newcommand{\be}{\begin{equation}}
\newcommand{\ee}{\end{equation}}
\newcommand{\beano}{\begin{eqnarray*}}
\newcommand{\eeano}{\end{eqnarray*}}
\newcommand{\inp}[2]{\langle {#1} ,\,{#2} \rangle}
\def\bmatrix#1{\left[ \begin{matrix} #1 \end{matrix} \right]}
\def \noin{\noindent}
\newcommand{\evenindex}{\Pi_e}



\def \R{{\mathbb R}}
\def \C{{\mathbb C}}
\def \K{{\mathbb K}}
\def \J{{\mathbb J}}
\def \Lb{\mathrm{L}}

\def \T{{\mathbb T}}
\def \Pb{\mathrm{P}}
\def \N{{\mathbb N}}
\def \Ib{\mathrm{I}}
\def \Ls{{\Lambda}_{m-1}}
\def \Gb{\mathrm{G}}
\def \Hb{\mathrm{H}}
\def \Lam{{\Lambda_{m}}}
\def \Qb{\mathrm{Q}}
\def \Rb{\mathrm{R}}
\def \Mb{\mathrm{M}}
\def \norm{\nrm{\cdot}\equiv \nrm{\cdot}}

\def \P{{\mathbb P}_m(\C^{n\times n})}
\def \A{{{\mathbb P}_1(\C^{n\times n})}}
\def \H{{\mathbb H}}
\def \L{{\mathbb L}}
\def \G{{\mathcal G}}
\def \S{{\mathbb S}}
\def \sigmin{\sigma_{\min}}
\def \elam{\sigma_{\epsilon}}
\def \slam{\sigma^{\S}_{\epsilon}}
\def \Ib{\mathrm{I}}
\def \Tb{\mathrm{T}}
\def \d{{\delta}}

\def \Lb{\mathrm{L}}
\def \N{{\mathbb N}}
\def \Ls{{\Lambda}_{m-1}}
\def \Gb{\mathrm{G}}
\def \Hb{\mathrm{H}}
\def \Delta{\triangle}
\def \Rar{\Rightarrow}
\def \p{{\mathsf{p}(\lam; v)}}

\def \D{{\mathbb D}}

\def \tr{\mathrm{Tr}}
\def \cond{\mathrm{cond}}
\def \lam{\lambda}
\def \sig{\sigma}
\def \sign{\mathrm{sign}}

\def \ep{\epsilon}
\def \diag{\mathrm{diag}}
\def \rev{\mathrm{rev}}
\def \vec{\mathrm{vec}}

\def \sk{\mathsf{skew}}
\def \sy{\mathsf{sym}}
\def \en{\mathrm{even}}
\def \odd{\mathrm{odd}}
\def \rank{\mathrm{rank}}
\def \pf{{\bf Proof: }}
\def \dist{\mathrm{dist}}
\def \rar{\rightarrow}

\def \rank{\mathrm{rank}}
\def \pf{{\bf Proof: }}
\def \dist{\mathrm{dist}}
\def \Re{\mathsf{Re}}
\def \Im{\mathsf{Im}}
\def \re{\mathsf{re}}
\def \im{\mathsf{im}}

\def \sym{\mathsf{sym}}
\def \sksym{\mathsf{skew\mbox{-}sym}}
\def \odd{\mathsf{odd}}
\def \even{\mathsf{even}}
\def \herm{\mathsf{Herm}}
\def \skherm{\mathsf{skew\mbox{-}Herm}}
\def \str{\mathrm{ Struct}}
\def \eproof{$\blacksquare$}
\def \proof{\noin\pf}

\def \bS{{\bf S}}
\def \cA{{\cal A}}
\def \E{{\mathcal E}}
\def \X{{\mathcal X}}
\def \F{{\mathcal F}}
\def \tr{\mathrm{Tr}}
\def \range{\mathrm{Range}}

\def \pal{\mathrm{palindromic}}
\def \palpen{\mathrm{palindromic~~ pencil}}
\def \palpoly{\mathrm{palindromic~~ polynomial}}
\def \hodd{H\mbox{-}\odd}
\def \heven{H\mbox{-}\even}


\title{On backward errors of structured polynomial eigenproblems solved by structure preserving linearizations}
\author{ Bibhas Adhikari\thanks{Department of Mathematics,
IIT Guwahati, India, E-mail:
bibhas.adhikari@gmail.com} \, and Rafikul Alam\thanks{
   Department of Mathematics,
IIT Guwahati, India, E-mail:
rafik@iitg.ernet.in, rafikul@yahoo.com,   Fax:
+91-361-2690762/2582649. }  }
\date{}

\maketitle \thispagestyle{empty}

{\small \noin{\bf Abstract.} First, we derive explicit computable expressions of
structured backward errors of approximate eigenelements of structured
matrix polynomials including symmetric, skew-symmetric,
Hermitian, skew-Hermitian, even and odd polynomials. We also determine minimal
structured perturbations for which  approximate eigenelements are
exact eigenelements of the perturbed polynomials. Next, we analyze the effect of structure preserving linearizations of structured matrix polynomials on the structured backward errors of approximate eigenelements. We identify structure preserving linearizations which have almost no adverse effect on the structured backward errors of approximate eigenelements of the polynomials. Finally, we analyze structured pseudospectra of a structured matrix polynomial and establish a partial equality between unstructured and structured pseudospectra. }

\vone \noin{\bf Keywords.} Structured matrix polynomial, structured backward
error, pseudospectrum, structured linearization.

\noin{\bf AMS subject classifications.} 65F15, 15A57, 15A18, 65F35

\section{Introduction} Consider a matrix polynomial $\Pb(z) :=\sum_{j=0}^m
z^j A_j$ of degree $m,$ where $ A_j\in \C^{n\times n} \,\,\mbox{and}\,\, A_m\neq
0.$ We assume that $\Pb$ is regular, that is, $\det(\Pb(z)) \neq 0$ for some $ z\in \C.$ We say that $\lam \in \C$ is an eigenvalue of $\Pb$ if $\det(\Pb(\lam)) =0.$ A nonzero vector $x \in \C^n$ (resp., $y \in \C^n$) that satisfies  $\Pb(\lam)x=0$ (resp., $y^H\Pb(\lam)=0$) is called a right (resp., left) eigenvector of $\Pb$ corresponding to the eigenvalue $\lam.$ The standard approach to computing eigenelements of $\Pb$ is to
convert $\Pb$ into an equivalent linear polynomial $\Lb$, called a
linearization of $\Pb,$  and employ a
numerically backward stable algorithm to compute the eigenelements
of $\Lb,$ where $\Lb(z) :=  zX + Y, \,  X\in  \C^{mn\times mn}$ and $  Y\in\C^{mn\times mn} .$
It is well known that a matrix polynomial admits several linearizations. In fact, it is shown in~\cite{mackey3, mackey:thesis} that potential linearizations of a matrix polynomial form a vector space. Thus choosing an optimal (in some sense) linearization of $\Pb$ is an important first step towards computing eigenelements of $\Pb.$  In general, a linearization of $\Pb$ can have an adverse effect on the conditioning of the eigenvalues of $\Pb$~(see, \cite{higcon}). Hence by analyzing the condition numbers of eigenvalues of linearizations, potential linearizations of $\Pb$ have been identified in \cite{higcon} whose eigenvalues are almost as sensitive to perturbations as that of $\Pb.$ Further, it is shown in \cite{higg1} that these linearizations  are consistent with the backward errors of approximate eigenelements in the sense that they nearly minimize the backward errors.

Polynomial eigenvalue problems that occur in many applications possess some distinctive structures (e.g., Hermitian, even, odd and palindromic) which in turn induce certain spectral symmetries on the eigenvalues of the matrix polynomials~(see, \cite{mackey2,schraoder:thesis, volker4,X, xu} and the references therein). With a view to preserving spectral symmetry in the computed eigenvalues (and possibly improved accuracy), there has been a lot of interests in developing structured preserving algorithms~(see, \cite{hillmm04,volker3,schraoder:thesis,mmmm08} and the references therein). Since linearization is the standard way to solve a polynomial eigenvalue problem, for a structured matrix polynomial it is therefore necessary to choose a structured linearization and then solve the linear problem by a backward stable structure preserving algorithm. For the accuracy assessment of computed solution, it is therefore important to understand the sensitivity of eigenvalues of a structured matrix polynomial with respect to structure preserving perturbations. Also it is equally important to know the structured backward errors of approximate eigenelements of a structured matrix polynomial. Moreover, for a variety of structured polynomials such as symmetric, skew-symmetric, Hermitian, skew-Hermitian, even, odd and palindromic polynomials, there are infinitely many structured linearizations, see~\cite{mackey1, mackey2}. This poses a genuine problem of choosing one linearization over the other. For computational purposes, it is highly desirable to know how different structured linearizations affect the accuracy of computed eigenelements. Thus the selection of an optimal or a near optimal structured linearization is an important step in the solution process of a structured polynomial eigenvalue problem. The sensitivity analysis of eigenvalues of structured matrix polynomials with respect to structure preserving perturbation has been investigated in~\cite{BaR09:cond}. It also provides a recipe for choosing  structured linearizations whose eigenvalues are almost as sensitive to structure preserving perturbations as that of the structured matrix polynomials.

To complete the investigation, in this paper we analyze structured backward errors of approximate eigenelements of symmetric, skew-symmetric, Hermitian, skew-Hermitian, $T$-even, $T$-odd, $H$-even and  $H$-odd polynomials. These structures are defined in Table~\ref{table:strpoleigenpair}. The main  contribution of this paper is as follows.

 First, we derive explicit computable expressions for the structured backward errors of approximate eigenelements of structured matrix polynomials. We also construct a  minimal structured perturbation so that an approximate eigenelement is the exact eigenelement of the structured perturbed polynomial. These results generalize similar results in \cite{BaR09:pencil} obtained for structured matrix pencils.

Second, we consider structured linearizations that preserve spectral symmetry of a structured matrix polynomial and compare the structured backward errors of approximate eigenelements with that of the structured polynomial. For example, a $T$-even matrix polynomial admits $T$-even as well as $T$-odd linearizations both of which preserve the spectral symmetry of the $T$-even polynomial. Based on these results we identify structured linearizations which are  optimal in the sense that the structured backward errors of  approximate eigenelements of the linearizations  are bounded above and below by a small constant multiple of that of the structured polynomials. We show that these linearizations are consistent with the choice of linearizations discussed in \cite{BaR09:cond} by analyzing structured condition numbers of eigenvalues.

Third, we show that the effect of structure preserving linearization on the structured backward errors of approximate eigenelements is almost harmless for a wide class of structured linearizations.  We show that bad effect, if any, of a structure preserving linearization can be neutralized by considering a complementary structured linearization.
For example, when $\Pb$ is a $T$-even polynomial, we show that any $T$-even linearization is optimal for eigenvalues $\lam$ of $\Pb$ such that $|\lam| \leq 1,$ and any $T$-odd linearization is optimal for eigenvalues $\lam$ such that $|\lam| \geq 1.$ In such  a case, we show that the backward error of an approximate eigenelement of the linearization  differ from that of $\Pb$ by no more than a factor of $2.$ We show that similar results hold for other structured polynomials as well. In contrast, it is shown in~\cite{BaR09:cond} that the condition numbers of eigenvalues of these optimal linearizations  differ from that of the polynomial by a factor of a  constant whose size could of the order of the degree of the matrix polynomial.

Finally, we analyze structured pseudospectra of structured matrix polynomials and establish a partial equality between structured and unstructured pseudospectra. Similar study for palindromic matrix polynomials has been carried out in~\cite{Ba09}, see also~\cite{Ba:thesis, sbora}.

The rest of the paper is organized as follows. In section~\ref{sec1}, we review structured polynomials and their spectral symmetries. In section~\ref{sec2}, we analyze structured backward errors of approximate eigenpairs of structured polynomials. In section~\ref{sec3},  we analyze the effect of structure preserving linearizations  on the backward errors of approximate eigenelements of structure polynomials and provide a recipe for choosing optimal  linearizations. Finally, in section~\ref{chapter4:pseudo}, we consider structured pseudospectra of structured matrix polynomials.

\section{Structured matrix polynomials }\label{sec1}
We consider matrix polynomial of degree $m$ of the form
$\Pb(z):=\sum_{j=0}^m z^jA_j,$ where $ A_j\in\C^{n \times n}$ and $A_m \neq 0.$ Let
$\P$ denote the vector space of matrix polynomials of degree at
most $m.$ The spectrum of a regular polynomial $\Pb\in\P,$ denoted by
$\sig(\Pb),$ is given by $\sig(\Pb) := \{ z\in\C : \det(\Pb(z)) =0
\}.$ Strictly, speaking $\sig(\Pb)$ consists of finite eigenvalues
of $\Pb.$ If the leading coefficient of $\Pb$ is singular then
$\Pb$ has an infinite eigenvalue. In this paper, we consider only
finite eigenvalues of matrix polynomials. An infinite eigenvalue
of $\Pb,$ if any, can easily be analyzed by considering the {\it
reverse} polynomial of $\Pb$~(see~\cite{sa2}).  We say that
$(\lam, x, y)$ is an eigentriple of $\Pb$ if $\lam $ is an
eigenvalue of $\Pb$ and, $x$ and $y$ are the corresponding nonzero
right and left eigenvectors, that is, $\Pb(\lam) x = 0$ and
$y^H\Pb(\lam) =0.$

We denote the transpose and conjugate transpose of a matrix $A$ by
$A^T$ and $A^H,$ respectively. Define the map $\P \rightarrow \P,
\, \Pb \mapsto \Pb^*$ given by $\Pb^*(z):=\sum_{j=0}^m z^j A_j^*, $
where $A^*=A^T$ or $A^*=A^H.$  The map $\Pb \mapsto \Pb^*$
can be used to define interesting structured matrix polynomials
such as symmetric, skew-symmetric, Hermitian, skew-Hermitian,
$*$-even and $*$-odd matrix polynomials. These structures are
defined in Table~\ref{table:strpoleigenpair}. The table also shows the eigentriples as well as the
spectral symmetries of the eigenvalues, see also~\cite{mackey2}. We denote the set of structured
polynomials having one of the structures given in
Table~\ref{table:strpoleigenpair} by $\S.$ By writing a
pair $(\lam, \mu)$ in the third column of
Table~\ref{table:strpoleigenpair} we mean that if $\lam$ is an
eigenvalue of $\Pb$ then so is $\mu.$ Notice that the eigenvalues of
Hermitian and skew-Hermitian polynomials have the same spectral
symmetry. Similarly, the eigenvalues of $*$-even and $*$-odd
polynomials have the same spectral symmetry, where $* \in \{T,
H\}.$

\begin{table}[h]
\begin{center}\renewcommand{\arraystretch}{1.2}
\begin{tabular}{|c|c|c|c|}
  \hline
  \textbf{$\S$} & \textbf{Condition} &\textbf{spectral symmetry} & \textbf{eigentriple}  \\
  \hline \hline
  symmetric  & $\Pb^{T}(z)= \Pb(z),\forall z\in \C$ & $\lam$ & $(\lambda,x,\overline{x})$ \\ \cline{1-2}
skew-symmetric & $\Pb^{T}(z)= -\Pb(z),\forall z\in \C$ & & \\
  \hline
  $T$-even  & $ \Pb^{T}(z)=\Pb(-z),\forall z\in \C $ & $(\lambda,-\lambda)$ & $(\lambda,x,\overline{y}),(-\lambda,y,\overline{x})$
  \\ \cline{1-2}
$T$-odd & $\Pb^{T}(z)=-\Pb(-z),\forall z\in \C$ & & \\
      \hline
    \hline
   Hermitian & $\Pb^{H}(z)=\Pb(z),\forall z\in \C$ & $(\lambda,\overline{\lambda})$ & $(\lambda,x,y), (\overline{\lambda},y,x)$
   \\\cline{1-2}
    skew-hermitian & $\Pb^{H}(z)=-\Pb(z),\forall z\in \C$ && \\
     \hline
    $H$-even  & $\Pb^{H}(z)=\Pb(-z),\forall z\in \C$ & $(\lambda,-\overline{\lambda})$ & $(\lambda,x,y), (-\overline{\lambda},y,x)$
    \\\cline{1-2}
    $H$-odd & $\Pb^{H}(z)=-\Pb(-z),\forall z\in \C$ && \\
       \hline
\end{tabular}
\caption{\label{table:strpoleigenpair} Spectral symmetries of
structured polynomials.}
\end{center}
\end{table}

Let $ \Pb \in \S$ be regular. With a view to obtaining structured backward error of
 $(\lam,x)\in\C\times\C^n$ with $x^Hx=1$ as an approximate eigenpair of $\Pb,$  we now show that there always exists
a polynomial $\Delta\Pb\in\S$ such that $(\lam,x)$ is a right
eigenpair of $\Pb+\Delta\Pb,$ that is,
$(\Pb(\lam)+\Delta\Pb(\lam)) x=0.$ Recall that $\S$ denotes the set of structured
polynomials having one of the structures given in
Table~\ref{table:strpoleigenpair}. In short, we write $ \S
\in \{ \sym, \sksym, \herm, \skherm, T\mbox{-}\even,
T\mbox{-}\odd, H\mbox{-}\even, H\mbox{-}\odd\}.$

\begin{theorem}\label{exist1} Let $ \S
\in \{ \sym, \sksym, \herm, \skherm, T\mbox{-}\even,
T\mbox{-}\odd, H\mbox{-}\even, H\mbox{-}\odd\}$ and  $ \Pb \in \S$
be given by $\Pb(z)=\sum_{j=0}^m z^j A_j.$ Let
$(\lambda,x)\in\C\times\C^n$ be such that $x^Hx=1.$ Set
$r=-\Pb(\lambda)x,$  $\Lam :=[1, \, \lam, \, \hdots, \, \lam^m]^T$ and $P_x := I-xx^H.$
 Define \beano \Delta A_j &:=& \left\{%
\begin{array}{ll}
    -\overline{x}x^{T}A_jxx^{H}+\frac{\overline{\lambda^{j}}}{\|\Lam\|_2^2}[\overline{x}r^{T} +
rx^{H} - 2 (r^Tx) \overline{x}x^H ], & \hbox{if~$A_j=A_j^T$,} \\\\
    -\frac{\overline{\lambda^{j}}}{\|\Lam\|_2^2}[\overline{x}r^{T} - rx^{H}], & \hbox{if~$A_j=-A_j^T$,} \\
\end{array}%
\right.\\
\Delta A_j &:=& \left\{%
\begin{array}{ll}
    -xx^{H}A_jxx^{H}+\frac{1}{\|\Lam\|_2^2}[\lambda^{j}xr^{H}P_x +
\overline{\lambda^{j}} P_x rx^{H}], & \hbox{if~$A_j=A_j^H$,} \\\\
    -xx^{H}A_jxx^{H}-\frac{1}{\|\Lam\|_2^2}[\lambda^{j}xr^{H}P_x -
\overline{\lambda^{j}} P_x rx^{H}], & \hbox{if~$A_j=-A_j^H$,} \\
\end{array}%
\right.\eeano and consider the polynomial $\Delta \Pb (z) :=
\sum_{j=0}^m z^j \Delta A_j.$ Then $\Pb (\lam)x + \Delta
\Pb(\lam)x =0$ and $\Delta \Pb \in \S.$
\end{theorem}

\noin \pf The proof is computational and is easy to
check.$\blacksquare$

\section{Structured backward errors}\label{sec2}
Backward errors of approximate eigenelements of  regular matrix polynomials have been systematically analyzed and computable expressions for the backward errors have been derived by Tisseur in~\cite{fran7} . For our purpose, we require a different norm setup for matrix polynomials.   We equip $\P$ with a norm so that the resulting normed
linear space can be used for perturbation analysis of matrix
polynomials.  Let $\Pb \in \P$ be given by $\Pb(z) :=
\sum^m_{j=0}A_j z^j.$  We define $$ \nrm{\Pb}_M :=
\left(\sum_{j=0}^m \|A_j\|_M^2\right)^{1/2},$$ where $\|A\|_M$ denotes the Frobenius norm when $M=F$ and the spectral norm when $M=2.$ Accordingly, we say that
$\nrm{\cdot}_F$ is the Frobenius norm and $\nrm{\cdot}_2$ is the spectral norm  on $\P$.
See~\cite{sa2, safique:thesis} for more on norms of matrix
polynomials.

Let $ (\lam, x) \in \C\times \C^n$ be such that $x^Hx=1$ and $\Pb \in \P$ be  regular.
We denote the backward error of $(\lam, x)$ as an approximate  eigenelement of $\Pb$ by $\eta_M(\lam, x, \Pb)$ given by $$\eta_M(\lam, x,\Pb) := \inf\limits_{\Delta
\Pb\in\P}\{ \nrm{\Delta\Pb}_M :
 \Pb(\lam)x +\Delta\Pb(\lam)x = 0\}.$$
Setting $r := -\Pb(\lam)x$ and $\Lam :=[1, \lam, \ldots,
\lam^m]^T,$ it is easily seen that \be\label{unstr:bckerr}
\eta_M(\lam, x, \Pb) = \frac{\|r\|_2}{\|x\|_2\|\Lam\|_2}\ee for $M=F$ as well as $M=2.$ Indeed, defining $\dm{\Delta{A_j} :=
\frac{\overline{\lam^j} rx^H}{x^Hx \|\Lam\|_2^2}},~j=0:m,$ and
considering the polynomial $\Delta\Pb(z) := \sum_{j=0}^m z^j \Delta
A_j,$ we have $\nrm{\Delta\Pb}_M = \|r\|_2/{\|x\|_2\|\Lam\|_2} $ and $\Pb(\lam)x+\Delta\Pb(\lam)x =0.$
Consequently,  for simplicity of notation, we denote $\eta_M(\lam,
x, \Pb)$   by $\eta(\lam, x,
\Pb).$

Now suppose that $\Pb\in\S.$ Then treating $(\lam, x) $ as an
approximate eigenelement of $\Pb,$  we define the structured backward
error of $(\lam, x)$ by $$\eta_M^{\S}(\lam, x, \Pb) :=
\inf\limits_{\Delta \Pb\in\S}\{ \nrm{\Delta\Pb}_M : \Pb(\lam)x
+\Delta\Pb(\lam)x = 0\}.$$  In
view of Theorem~\ref{exist1}, it follows that
 $\eta(\lam, x, \Pb) \leq \eta^{\S}_M(\lam,
x,\Pb) < \infty.$  Structured backward errors of approximate eigenelements of structured matrix pencils have been systematically analyzed and computable expressions of the structured backward errors have been derived in~\cite{BaR09:pencil}. In this section we generalize these results to the case of structured matrix polynomials.

As we shall see, determining $\eta^\S_2(\lam, x, \Pb)$ is much more difficult than determining $\eta_F^\S(\lam, x,
\Pb)$ and requires  solution of norm preserving dilation problem
for matrices. The Davis-Kahan-Weinberger solutions of norm
preserving dilation problem given below will play an important
role in the subsequent development. Let $A,B,C$ and $D$ be
matrices of appropriate sizes. Then the following result holds.

\begin{theorem}[Davis-Kahan-Weinberger, \cite{D}] \label{dkw}
Let $A, B, C$ satisfy $ \left\|\bmatrix{ A\\
B}\right\|_2 = \mu$ and \\ $ \left\|\bmatrix{ A & C}\right\|_2
= \mu.$ Then   there exists $D$ such that $\left\|\bmatrix{  A & C \\
  B & D
}\right\|_2 = \mu.$ Indeed, those $D$ which have this property are
exactly those of the form $$D= - KA^HL +
\mu(I-KK^H)^{1/2}Z(I-L^HL)^{1/2},$$ where $K^H :=(\mu^2 \Ib -A^HA
)^{-1/2}B^H,~ L :=(\mu^2 \Ib -AA^H)^{-1/2}C$ and $Z$ is an
arbitrary contraction, that is, $\|Z\|_2 \leq 1.$  $\blacksquare$
\end{theorem}

For a more general version of the above result, see~\cite{D}.

\subsection{Symmetric and skew-symmetric polynomials}
We now derive structured backward error of $(\lam, x) \in
\C\times \C^n$ as an approximate eigenpair of symmetric and
skew-symmetric matrix polynomials. We also derive minimal
structured perturbations so that $(\lam, x)$ is an exact eigenpair of
the perturbed polynomials. First, we consider symmetric matrix polynomials. Note that a
matrix polynomial $\Pb \in \P$ is symmetric if and only if all the
coefficient matrices of $\Pb$ are symmetric. For a symmetric
matrix polynomial, we have the following result.

\begin{theorem}\label{tsym} Let $\S$ denote the set of symmetric matrix
polynomials in $\P$ and let $\Pb\in\S.$ Let $(\lam, x) \in
\C\times \C^n$ be such that $x^Hx =1.$ Set $r:=-\Pb (\lambda)x,$ $P_x
:= I-xx^H$ and $\Lam :=[1, \lam, \ldots, \lam_m]^T.$ Then we have
$$\eta_F^{\S}(\lambda,x,\Pb)=\frac{\sqrt{2\|r\|_{2}^{2}-|x^{T}r|^{2}}}
{\|\Lam\|_2}\leq \sqrt{2} \eta(\lam,x,\Pb) \,\mbox{ and } \,
\eta_2^{\S}(\lambda,x,\Pb)=\eta(\lam,x,\Pb).$$ Set $\Delta A_j
:=\frac{\overline{\lambda^{j}}}{\|\Lam\|_2^2}[\overline{x}r^T+rx^H-(r^Tx)\overline{x}x^H],
\, \,  \, j = 0:m,$ and consider the polynomial $\Delta \Pb(z):=
\sum_{j=0}^m z^j \Delta A_j.$ Then $\Delta \Pb $ is a unique
polynomial such that $\Delta\Pb\in\S, \,\,\Delta
\Pb(\lam)x+\Pb(\lam)x=0$ and $\nrm{\Delta
\Pb}_F=\eta_F^{\S}(\lambda,x,\Pb).$ Further, define
$$\Delta A_j :=\frac{\overline{\lambda^{j}}}{\|\Lam\|_2^2}[\overline{x}r^T+rx^H-(r^Tx)\overline{x}x^H]-
\frac{\overline{\lam^j}~\overline{x^Tr}~P_x^T rr^T
P_x}{\|\Lam\|_2^2~(\|r\|_2^2 - |x^Tr|^2)}$$ and consider the
polynomial $\Delta \Pb(z) := \sum_{j=0}^m z^j \Delta A_j.$ Then
$\Delta \Pb \in\S, \, \Delta \Pb(\lam)x+\Pb(\lam)x=0$ and
$\nrm{\Delta \Pb}_2=\eta_2^{\S}(\lambda,x,\Pb).$
\end{theorem}

\noin\pf In view of Theorem~\ref{exist1}, let $\Delta \Pb\in \S$ given by $\Delta \Pb(z)
:= \sum^m_{j=0}\Delta A_jz^j$  be such that $\Pb(\lam)x+\Delta \Pb(\lam)x=0.$  Let $Q_{1}\in\C^{n\times (n-1)}$   be
such that the matrix $Q=[x~~Q_{1}]$ is unitary.
Then $\widetilde{\Delta A_j} :=Q^{T}\Delta A_jQ=\left(%
\begin{array}{cc}
  a_{jj} & a_{j}^{T} \\
  a_{j} & X_{j} \\
\end{array}%
\right),$ where $X_j=X_j^T$ is of size $n-1.$ Since
$\overline{Q}Q^T=I,$ we have
$$\overline{Q}(\Delta \Pb (\lam))Q^H x=r\Rightarrow
(\Delta \Pb (\lam)) Q^H x=Q^Tr=\left(%
\begin{array}{c}
  x^{T}r \\
  Q^{T}_{1}r \\
\end{array}%
\right)$$ As $Q^H x=e_1,$ the first column of the identity matrix,
we have $\left(%
\begin{array}{c}
  \sum_{j=0}^m\lambda^{j}a_{jj} \\
  \sum_{j=0}^m\lambda^{j}a_{j} \\
\end{array}%
\right)=\left(%
\begin{array}{c}
  x^{T}r \\
  Q^{T}_{1}r \\
\end{array}%
\right).$ Hence the minimum norm solutions are
$a_{j}=\frac{\overline{\lambda^{j}}Q^{T}_{1}r}{\|\Lam\|_2^2}$ and
$
a_{jj}=\frac{\overline{\lambda^{j}}x^{T}r}{\|\Lam\|_2^2},~~j=0:m.$
Consequently, we have
\be\label{bck:tsympoly1}\widetilde{\Delta A}_j=\left(%
\begin{array}{cc}
  \frac{\overline{\lambda^{j}}x^{T}r}{\|\Lam\|_2^2} &
  (\frac{\overline{\lambda^{j}}Q^{T}_{1}r}{\|\Lam\|_2^2})^{T} \\[8pt]
  \frac{\overline{\lambda^{j}}Q^{T}_{1}r}{\|\Lam\|_2^2} & X_{j} \\
\end{array}%
\right).\ee This shows that the Frobenius norm of
$\widetilde{\Delta A_j}$ is minimized when $X_j=0.$ Hence we have
$\|\Delta A_j\|_F^2=\|\widetilde{\Delta
A_j}\|_F^2=|a_{jj}|^{2}+2\|a_{j}\|_{2}^{2}.$ Since $
\overline{Q}_{1}Q_{1}^{T}=I-\overline{x}x^{T},$ we have
$$\eta_F^\S(\lam,x,\Pb) = \sqrt{\frac{|x^{T}r|^{2}}{\|\Lam\|_2^2}+
\frac{2\|(I-\overline{x}x^{T})r\|_{2}^{2}}{\|\Lam\|_2^2}}=
\frac{\sqrt{2\|r\|_{2}^{2}-|x^{T}r|^{2}}} {\|\Lam\|_2}.$$
Now from (\ref{bck:tsympoly1}), we have $$ \Delta
A_j = [\overline{x} \,\,\, \overline{Q}_1]\left(%
\begin{array}{cc}
  \frac{\overline{\lambda^{j}}x^{T}r}{\|\Lam\|_2^2} &
  (\frac{\overline{\lambda^{j}}Q^{T}_{1}r}{\|\Lam\|_2^2})^{T} \\[8pt]
  \frac{\overline{\lambda^{j}}Q^{T}_{1}r}{\|\Lam\|_2^2} & 0 \\
\end{array}%
\right)\left(%
\begin{array}{c}
  x^H \\
  Q_1^H \\
\end{array}%
\right) \\
=
\frac{\overline{\lambda^{j}}}{\|\Lam\|_2^2}[\overline{x}r^T+rx^H-(r^Tx)\overline{x}x^H]$$
which gives the desired polynomial $\Delta \Pb$ for the Frobenius
norm.

For the spectral norm, we employ dilation result in
Theorem~\ref{dkw} to the matrix in (\ref{bck:tsympoly1}). Indeed,
for $ \mu_j := \frac{|\lam^j|~\|r\|_2}{\|\Lam\|_2^2},$ by
Theorem~\ref{dkw}, we have \beano  X_{j} =
-\frac{\overline{\lam^j}~\overline{x^Tr}~Q_1^Tr
(Q_1^Tr)^T}{\|\Lam\|_2^2~(\|r\|_2^2 - |x^Tr|^2)}, \, \, j=0:m,
\eeano which gives $\dm{
\eta_2^\S(\lam,x,\Pb)=\frac{\|r\|_2}{\|\Lam\|_2} =
\eta(\lam,x,\Pb)}.$  Putting $X_j$ in (\ref{bck:tsympoly1}) and after simplification we
have  $$\Delta
A_j=\frac{\overline{\lambda^{j}}}{\|\Lam\|_2^2}[\overline{x}r^T+rx^H-(r^Tx)\overline{x}x^H]-
\frac{\overline{\lam^j}~\overline{x^Tr}~P_x^T rr^T
P_x}{\|\Lam\|_2^2~(\|r\|_2^2 - |x^Tr|^2)}$$ which gives the
desired polynomial $\Delta \Pb$ for the spectral norm.
$\blacksquare$

\begin{remark}
If $|x^Tr|=\|r\|_2,$ then $\|Q_1^Tr\|_2=0$. Hence considering
$X_j=0, \, j=0:m,$ in the above proof we obtain the desired
results for the spectral norm. Note that in such a case we have
$\eta^\S_F(\lam, x, \Pb) = \sqrt{ 2} \, \eta(\lam, x, \Pb).$
\end{remark}

Observe that if $Y$ is symmetric and $Yx =0$ then  $ Y = P_x^T Z
P_x$ for some symmetric matrix $Z.$  Consequently, from the proof
Theorem~\ref{tsym},  we have $\overline{Q_j}X_jQ_j^H =P_x^T Z_j
P_x, j=0:m,$ for some symmetric matrices $Z_j.$ Hence  we have
following.

\begin{corollary}\label{gpol}
Let $\Pb$ be a symmetric matrix polynomial. For
$(\lambda,x)\in\C\times\C^{n}$ with $x^Hx=1,$ set $r
:=-\Pb(\lambda)x.$ Then there is a symmetric matrix polynomial
$\Qb$ such that $\Pb(\lam)x + \Qb(\lam)x=0$ if and only if
$\Qb(z)=\Delta\Pb(z) + P_x^T\Rb(z)P_x$ for some symmetric
polynomial $\Rb$, where $\Delta \Pb$ is the symmetric polynomial
given by $\Delta\Pb(z) :=\sum_{j=0}^m z^j \Delta A_j$ and
$$ \Delta A_j :=
\frac{\overline{\lambda^{j}}}{\|\Lam\|_2^2}[\overline{x}r^T+rx^H-(r^Tx)\overline{x}x^H],
\, j=0:m.$$
\end{corollary}

Next, we consider skew-symmetric matrix polynomials. Note that a
matrix polynomial $\Pb \in \P$ is skew-symmetric if and only if
all the coefficient matrices of $\Pb$ are skew-symmetric.  For
skew-symmetric matrix polynomials we have the following result.

\begin{theorem}\label{tskewsym}
Let $\S$ denote the set of skew-symmetric matrix polynomials in
$\P$  and let $\Pb \in\S.$  For $(\lambda,x)\in\C\times\C^n$ with
$x^Hx =1,$  set $r :=-\Pb(\lambda)x.$ Then  we have
$$\eta_F^{\S}(\lambda,x,\Pb)=\sqrt{2}~\eta(\lambda,x,\Pb),\,\,\eta_2^{\S}(\lambda,x,\Pb)= \eta(\lam,x,\Pb).$$
For the skew-symmetric polynomial $\Delta \Pb $ given in
Theorem~\ref{exist1}, we have $\Pb(\lambda) x+\Delta
\Pb(\lambda)x=0, \, \nrm {\Delta
\Pb}_F=\eta_F^{\S}(\lambda,x,\Pb)$ and $\nrm {\Delta
\Pb}_2=\eta_2^{\S}(\lambda,x,\Pb).$
\end{theorem}

\noin \pf The proof is the same as that of Theorem~\ref{tsym} except that  $\Delta A_j$ is
skew-symmetric for $ j =0:m.$ This gives
\be\label{bck:tsksympoly1} \widetilde{\Delta A_j} = \left(%
\begin{array}{cc}
  0 &
  -(\frac{\overline{\lambda^{j}}Q^{T}_{1}r}{\|\Lam\|_2^2})^{T} \\
  \frac{\overline{\lambda^{j}}Q^{T}_{1}r}{\|\Lam\|_2^2} & X_{j} \\
\end{array}%
\right).\ee Setting $X_j=0,$ we obtain the results for the Frobenius norm.

Setting $ \mu_j := \frac{|\lam^j|~\|r\|_2}{\|\Lam\|_2^2}$ and invoking Theorem~\ref{dkw}, it is easily seen that the spectral norm of
$\widetilde{\Delta A_j}$ in (\ref{bck:tsksympoly1}) is minimized when $X_j=0.$
Hence the desired results follow for the spectral norm. $\blacksquare$

\vone Note that if $Y$ is a skew-symmetric matrix and $Yx=0$ then
$Y = P_x^T ZP_x$ for some skew-symmetric matrix $Z.$ Hence we have
the following result.

\begin{corollary}\label{cor:tsksym}
Let $\Pb\in\P$ be a skew-symmetric matrix polynomial. For
$(\lambda,x)\in\C\times\C^{n}$ with $x^Hx=1,$ set $r
:=-\Pb(\lambda)x.$ Then there is a skew-symmetric matrix
polynomial $\Qb$ such that $\Pb(\lam)x + \Qb(\lam)x=0$ if and only
if $\Qb(z)=\Delta\Pb(z) + P_x^T\Rb(z)P_x$ for some skew-symmetric
polynomial $\Rb$, where $\Delta\Pb$ is the skew-symmetric
polynomial given by $\Delta\Pb(z) :=\sum_{j=0}^m z^i\Delta A_j$
and
$$\Delta A_j =   -\frac{\overline{\lambda^{j}}}{\|\Lam\|_2^2}[\overline{x}r^{T} - rx^{H}],\,\, j =0:m.$$
\end{corollary}

\subsection{T-even and T-odd matrix polynomials}
 For backward perturbation analysis of $T$-even and $T$-odd polynomials, we need the even
 index projection $\Pi_e : \C^{m+1} \rightarrow \C^{m+1}$ given by
\[
    \Pi_e ([x_0, \, x_1, \, \hdots, \, x_{m-1}, \, x_m]^T) :=
 \left\{%
 \begin{array}{ll}
     [ x_0,\, 0,\,  x_2,\,0,\, \ldots,\,\, x_{m-2},\, 0,\, x_m]^T, & \hbox{if $m$ is even}, \\{}
     [ x_0,\, 0,\, x_2,\,0,\, \ldots,\, 0,\, x_{m-1},\, 0]^T, & \hbox{if $m$ is odd.} \\
 \end{array} \right.
\] Note that $``0"$ is considered as even number. Observe that  $I - \Pi_e$ is the odd index
projection.

Recall that a matrix polynomial $\Pb \in \P$ given by $\Pb(z) :=
\sum^m_{j=0} A_j z^j$ is $T$-even  if and only if $A_j$ is
symmetric when $j$ is even (including $j=0$) and $A_j$ is
skew-symmetric when $j$ is odd. We have the following result for
$T$-even matrix polynomials.

\begin{theorem}\label{teven}
Let $\S$ denote the set of  $T$-even matrix polynomials in $\P.$
Let $\Pb\in\S$ and  $ (\lam, x) \in \C\times \C^n$ be such that
$x^Hx=1.$ Set $r :=-\Pb(\lambda)x,$ $P_x := I-xx^H$ and $ \Lam
:=[1, \lam, \ldots, \lam^m]^T.$ Then  we have
$$\eta_F^\S(\lambda,x,\Pb)=\sqrt{\frac{|x^{T}r|^{2}}{\|\Pi_e(\Lam)\|_2^2}+
2\frac{\|r\|_{2}^{2}-|x^Tr|^2}{\|\Lam\|_2^2}}, \,\,
\eta_2^\S(\lambda,x,\Pb) = \sqrt{\frac{|x^{T}r|^{2}}{\|\Pi_e
(\Lam)\|_2^2}+ \frac{\|r\|_{2}^{2}-|x^Tr|^2}{\|\Lam\|_2^2}}.$$

In particular, if $m$ is odd and $|\lam|=1$ then we have
$\eta_F^\S(\lambda,x,\Pb)= \sqrt{2}~\eta(\lam,x,\Pb)$ and $
\eta_2^\S(\lambda,x,\Pb) = \eta(\lam,x,\Pb).$

For $j =0 : m,$ define  $$ \dm{E_j := \left\{%
\begin{array}{ll}
    \dfrac{\overline{\lambda^{j}}}{\|\Pi_e(\Lam)\|_2^2} (x^{T}r)
\overline{x}x^H +\dfrac{\overline{\lam^{j}}}{\|\Lam\|_2^2}
[\overline{x} r^{T}P_x + P_x^T rx^{H}], & \hbox{ if $j$ is even, }
\\[8pt]
    \dfrac{\overline{\lam^{j}}}{\|\Lam\|_2^2}[P_x^T
rx^{H}-\overline{x} r^{T} P_x ], & \hbox{if $j$ is odd .} \\
\end{array}%
\right.}$$  Then $\Delta \Pb(z) :=\sum_{j=0}^m z^j E_j$ is a
unique $T$-even  polynomial in $\S$ such that $\Pb
(\lambda)x+\Delta \Pb (\lambda)x=0$ and $\nrm{\Delta\Pb}_F =
\eta^\S_F(\lam,x,\Pb).$ Further, for $j=0:m,$ defining
$$\Delta A_j := \left\{%
\begin{array}{ll}
    E_{j}
-\dfrac{\overline{\lam^{j}}~\overline{x^Tr}~P_x^T rr^T
P_x}{\|\Pi_e(\Lam)\|_2^2~(\|r\|_{2}^{2}-|x^Tr|^2)}, & \hbox{if $j$ is even, } \\
    E_j, & \hbox{if $j$ is odd, } \\
\end{array}%
\right. $$ we obtain a T-even polynomial $\Delta \Pb(z)
:=\sum_{j=0}^m z^j\Delta A_j$ in $\S$ such that $\Pb
(\lambda)x+\Delta \Pb (\lambda)x=0$ and $\nrm{\Delta\Pb}_2 =
\eta^\S_2(\lam,x,\Pb).$

\end{theorem}

\noin \pf  In view of Theorem~\ref{exist1}, let $\Delta \Pb\in \S$ be
such that $\Pb(\lam)x+\Delta \Pb(\lam)x=0.$ Assuming that $\Delta \Pb$ is given by $\Delta \Pb(z)
:= \sum^m_{j=0}\Delta A_jz^j,$ and arguing similarly as in the proofs of Theorems~\ref{tsym} and \ref{tskewsym}, we have
 $\widetilde{\Delta A_j}
=   \left(%
\begin{array}{cc}
  a_{jj} & a_{j}^{T} \\
  a_{j} & X_{j} \\
\end{array}%
\right),$ \, $X_{j}^T=X_{j}$ when $j$ is even, and $\widetilde{\Delta A_j}
=\left(%
\begin{array}{cc}
  0 & b_{j}^{T} \\
  -b_{j} & Y_{j} \\
\end{array}%
\right),~~Y_{j}^T=-Y_{j}$ when $j$ is odd. Consequently, we
have
$$\left(%
\begin{array}{c}
  \sum_{j} \lam^{j} a_{jj} \\
  \sum_{j\mbox{-even}} \lam^{j} a_{j}-\sum_{j\mbox{-odd}} \lam^{j} b_{j} \\
\end{array}%
\right)=\left(%
\begin{array}{c}
  x^{T}r \\
  Q^{T}_{1}r \\
\end{array}%
\right).$$ Hence the smallest norm solutions are
$a_{jj}=\frac{\overline{\lam^{j}}} {\|\Pi_e(\Lam)\|_2^2}~x^Tr,\,
a_{j}= \frac{\overline{\lam^{j}}}{\|\Lam\|_2^2}Q^{T}_1 r,\, b_{j}
=-\frac{\overline{\lam^{j}}}{\|\Lam\|_2^2} Q^{T}_1 r. $ Therefore,
we have \be\label{bck:tevenpoly1}\widetilde{\Delta A_j} = Q^T \Delta A_j Q = \left\{%
\begin{array}{ll}
     \left(%
\begin{array}{cc}
  \frac{\overline{\lam^{j}}}{\|\Pi_e(\Lam)\|_2^2}x^{T}r &
  (\frac{\overline{\lambda^{j}}Q^{T}_{1}r}{\|\Lam\|_2^2})^{T} \\\\
  \frac{\overline{\lambda^{j}}Q^{T}_{1}r}{\|\Lam\|_2^2} & X_{j} \\
\end{array}%
\right), & \hbox{if $j$ is even} \\[10pt] \\
   \left(%
\begin{array}{cc}
  0 &
  -(\frac{\overline{\lambda^{j}}Q^{T}_{1}r}{\|\Lam\|_2^2})^{T}
  \\\\
  \frac{\overline{\lambda^{j}}Q^{T}_{1}r}{\|\Lam\|_2^2} & Y_{j} \\
\end{array}%
\right), & \hbox{if $j$ is odd.} \\
\end{array}%
\right. \ee Setting $X_{j}=0=Y_{j}$ and using the fact that
$\overline{Q}_1Q_1^T=\Ib-\overline{x}x^T,$ we obtain the desired unique
$T$-even polynomial $\Delta \Pb(z) := \sum_{j=0}^m z^j E_j$
such that
$$\nrm{\Delta\Pb}_F = \eta_F^\S(\lambda,x,\Pb) = \sqrt{\frac{|x^{T}r|^{2}}{\|\Pi_e(\Lam)\|_2^2}+
2\frac{\|r\|_{2}^{2}-|x^Tr|^2}{\|\Lam\|_2^2}}.$$ When $m$ is odd
and $|\lam|=1,$ it is easily seen that $\|\Pi_e(\Lam)\|_2^2 =
\frac{1}{2}\|\Lam\|_2^2.$ Hence we have
$\eta_F^\S(\lam,x,\Pb)=\sqrt{2}\,\eta(\lam,x,\Pb).$

For the spectral norm, setting  $\mu_j :=
\sqrt{\frac{|\lambda^{j}|^2~|x^{T}r|^2}{\|\Pi_e(\Lam)\|_2^4}+
\frac{|\lam^{j}|^2~(\|r\|_{2}^{2}-|x^Tr|^2)}{\|\Lam\|_2^4}}$ when
$j$ is even, and $\mu_j :=\sqrt{
\frac{|\lam^{j}|^2~(\|r\|_{2}^{2}-|x^Tr|^2)}{\|\Lam\|_2^4}}$ when
$j$ is odd, and applying Theorem~\ref{dkw} to the matrices in
(\ref{bck:tevenpoly1}), we have  \beano X_{j} &=&
-\frac{\overline{\lam^{j}}~\overline{x^Tr}~Q_1^Tr
(Q_1^Tr)^T}{\|\Pi_e(\Lam)\|_2^2~(\|r\|_{2}^{2}-|x^Tr|^2)}\,\,
\mbox{and}\,\, Y_{j}=0.\eeano Consequently, we have  $\dm{
\eta_2^\S(\lambda,x,\Pb) = \sqrt{\frac{|x^{T}r|^{2}}{\|\Pi_e
(\Lam)\|_2^2}+ \frac{\|r\|_{2}^{2}-|x^Tr|^2}{\|\Lam\|_2^2}}
.}$ From~(\ref{bck:tevenpoly1}), we have \beano \Delta A_j = \left\{%
\begin{array}{ll}
     \dfrac{\overline{\lam^j} x^Tr \overline{x}x^H}{\|\Pi_e(\Lam)\|_2^2} + \dfrac{\overline{\lam^j}}{\|\lam\|_2^2}[\overline{x}r^T P_x +
P_x^Trx^H] + \overline{Q}_1X_jQ_1^H, & \hbox{if $j$ is even} \\[10pt]
   \dfrac{\overline{\lam^j}}{\|\lam\|_2^2}[P_x^Trx^H - \overline{x}r^T P_x  ] + \overline{Q}_1Y_jQ_1^{H}, & \hbox{if $j$ is odd.} \\
\end{array}%
\right.\eeano Substituting $X_j$ and $Y_j$ in $\Delta A_j$ we
obtain the desired $T$-even matrix polynomial $\Delta \Pb $ for the spectral norm. $\blacksquare$

\begin{remark}
If $|x^Tr|=\|r\|_2$ then $\|Q_1^Tr\|_2=0$. Hence considering
$X_{j}=0=Y_j$ in the above proof, we obtain the desired result for
the spectral norm. Note that in such a case we have
$\eta^\S_F(\lam, x, \Pb) = \sqrt{2} \, \eta^\S_2(\lam, x, \Pb) =
\sqrt{2} \eta(\lam, x, \Pb).$
\end{remark}

 Recall that when $A$ is symmetric (resp., skew-symmetric) and
$Ax=0$ then $A = P_x^T ZP_x$ for some symmetric (resp.,
skew-symmetric) matrix $Z.$ Consequently, from the proof
Theorem~\ref{teven} it follows that $ \Delta A_{j} := E_{j} +
P_x^T Z_{j} P_x, $ where $Z_{j} =Z_j^T$ when $j$ is even, and $
Z_{j}^T = -Z_{j}$ when $j$ is odd. Hence we have the following
result.

\begin{corollary} \label{cor:teven} Let $\Pb$ be a $T$-even matrix
polynomial in $\P.$ Let $( \lam, x ) \in \C\times \C^n$ be such
that $x^Hx=1.$ Then there is a $T$-even matrix polynomial $\Qb$
such that $\Pb(\lam)x +\Qb(\lam) x =0$ if and only if $\Qb(z) =
\Delta \Pb(z) + P_x^T \Rb(z) P_x$ for some $T$-even matrix
polynomial $\Rb \in \P,$ where $\Delta \Pb (z) := \sum_{j=0}^m E_j
z^j$ and $E_j$'s are given in Theorem~\ref{teven}.

\end{corollary}

Next, we consider backward error of $T$-odd polynomials. Observe
that a matrix polynomial $\Pb \in \P$ given by $\Pb(z) :=
\sum^m_{j=0} A_j z^j$ is $T$-odd  if and only if $A_j$ is
skew-symmetric when $j$ is even (including $j=0$) and $A_j$ is
symmetric when $j$ is odd.

\begin{theorem}\label{todd} Let $\S$ denote the set of  $T$-odd matrix polynomials in
$\P.$ Let  $\Pb\in\S$  and $ (\lam, x) \in \C\times \C^n$ be such
that $x^Hx=1.$ Set $r :=-\Pb(\lambda)x,$ $P_x := I-xx^H$ and $
\Lam :=[1, \lam, \ldots, \lam^m]^T.$ Then we have
\beano \eta_F^\S(\lambda,x,\Pb) &=& \left\{%
\begin{array}{ll}
    \sqrt{\frac{|x^{T}r|^{2}}{\|(I-\Pi_e)(\Lam)\|_2^2}+
2\frac{\|r\|_{2}^{2}-|x^Tr|^2}{\|\Lam\|_2^2}}, & \hbox{if $\lam\neq 0,$} \\
    \sqrt{2} \, \eta(\lam,x,\Pb), & \hbox{if $\lam=0,$} \\
\end{array}%
\right.\\ \eta_2^\S(\lambda,x,\Pb) &=& \left\{%
\begin{array}{ll}
    \sqrt{\frac{|x^{T}r|^{2}}{\|(I-\Pi_e)(\Lam)\|_2^2}+
\frac{\|r\|_{2}^{2}-|x^Tr|^2}{\|\Lam\|_2^2}}, & \hbox{if $\lam\neq 0,$} \\
   \eta(\lam,x,\Pb), & \hbox{if $\lam=0.$} \\
    \end{array}%
\right. \eeano

In particular, if  $m$ is odd and $|\lam|=1$ we have
$\eta_F^\S(\lambda,x,\Pb)= \sqrt{2}~\eta(\lam,x,\Pb)$ and $
\eta_2^\S(\lambda,x,\Pb) = \eta(\lam,x,\Pb).$ For $j =0:m, $ define  $$\dm{F_j := \left\{%
\begin{array}{ll}
    \dfrac{\overline{\lam^{j}}}{\|\Lam\|_2^2}[P_x^T
rx^H- \overline{x}r^T P_x ], & \hbox{if $j$ is even} \\[8pt]
    \dfrac{\overline{\lam^{j}}\overline{x}x^Trx^H}{\|(I-\Pi_e)(\Lam)\|_2^2}+\dfrac{\overline{\lam^{j}}}
{\|\Lam\|_2^2}[\overline{x}r^T P_x + P_x^T rx^H], & \hbox{if $j$ is odd.} \\
\end{array}%
\right.}$$ Then $\Delta \Pb(z) :=\sum_{j=0}^m z^j F_j$ is a unique
$T$-odd polynomial in $\S$ such that $\Pb (\lambda)x+\Delta \Pb
(\lambda)x=0$ and $\nrm{\Delta\Pb}_F = \eta^\S_F(\lam,x,\Pb).$

Further, for $j =0:m, $ define $ \Delta A_{j} := F_{j}$ when  $j$
is even, and $$\Delta A_{j} := F_{j} -\frac{\overline{\lam^{j}} \,
\overline{x^Tr} P_x^T rr^T
P_x}{\|(I-\Pi_e)(\Lam)\|_2^2(\|r\|_2^2-|x^Tr|^2)}$$ when $j$ is odd.
Then  $\Delta \Pb(z) :=\sum_{j=0}^m z^j\Delta A_j$ is a $T$-odd
polynomial in $\S$ such that $\Pb (\lambda)x+\Delta \Pb
(\lambda)x=0$ and $\nrm{\Delta\Pb}_2 = \eta^\S_2(\lam,x,\Pb).$
\end{theorem}

\noin \pf The desired results follow from the proof of
Theorem~\ref{teven} by interchanging the role of $\Delta A_{j}$
for even $j$ and odd $j.$ $\blacksquare$

We have the following results whose proof is immediate.

\begin{corollary} \label{cor:teven} Let $\Pb$ be a $T$-odd matrix
polynomial in $\P.$ Let $( \lam, x ) \in \C\times \C^n$ be such
that  $x^Hx=1.$ Then there is a $T$-odd matrix polynomial $\Qb$
such that $\Pb(\lam)x +\Qb(\lam) x =0$ if and only if $\Qb(z) =
\Delta \Pb(z) + P_x^T \Rb(z) P_x$ for some $T$-odd matrix
polynomial $\Rb \in \P,$ where $\Delta \Pb (z) := \sum_{j=0}^m F_j
z^j$ and $F_j$'s are given in Theorem~\ref{todd}.
\end{corollary}

\subsection{Hermitian and skew-Hermitian matrix polynomials}
We now consider structured  backward errors of approximate eigenelements  of
Hermitian and skew-Hermitian matrix polynomials. We proceed as
follows. Let $\S \subset \P$ and $\omega \in \C$ be such that
$|\omega|=1.$ We set $ \S_{\omega} := \{ \omega \Pb : \Pb \in \S\}.$
Then for $\Pb \in \P,$ it is easily seen that \be \label{eq:rotb}
\eta^{\S}_F(\lam, x, \Pb)= \eta^{\S_{\omega}}_F(\lam, x, \omega \Pb)
\mbox{ and } \eta^{\S}_2(\lam, x, \Pb)= \eta^{\S_{\omega}}_2(\lam, x,
\omega \Pb). \ee

Note that a matrix polynomial $\Pb \in \P$ is Hermitian (resp.,
skew-Hermitian) if and only if all the coefficient matrices of
$\Pb$ are Hermitian (resp., skew-Hermitian).  Let $\herm$ and
$\skherm,$ respectively, denote the set of Hermitian and
skew-Hermitian matrix polynomials in $\P.$ Then noting that a
matrix $X \in \C^{n\times n}$ is Hermitian if and only if  $ i X$
is skew-Hermitian, it easily seen that  the maps
\be\label{eq:rotp} \herm \longrightarrow \skherm,
 \, \Pb \longmapsto i \Pb\,\, \mbox{ and } \,\,  \skherm \longrightarrow \herm,
 \, \Qb \longmapsto i \Qb \ee are isometric isomorphisms. Thus, in view of (\ref{eq:rotb})
and (\ref{eq:rotp}), it follows that the structured backward error
of $(\lam, x)$ as an approximate eigenpair of a skew-Hermitian
polynomial can be obtained from the structured backward error of
$(\lam, x)$ as an approximate eigenpair of a Hermitian matrix
polynomial and vice-versa. We therefore analyze structured
backward perturbation of Hermitian matrix polynomials.

For $x \in \C^n,$ we denote by $\Re(x)$ and $\Im(x),$
respectively,  the real and the imaginary parts of $x.$ Then we
have $x = \Re(x) +i \Im(x).$ We denote the real and imaginary part
of a complex number $z \in \C$ by $\re(z)$ and $\im(z),$
respectively. We denote the Moore-Penrose pseudo-inverse of $A$ by $A^\dag$ and the canonical basis of $\C^{m+1}$
by $e_j, j=0:m.$

\begin{theorem}\label{hherm}  Let $\herm$ denote the set of Hermitian
matrix polynomials in $\P.$  Let  $\Pb\in\herm$  and $ (\lam, x)
\in \C\times \C^n$ be such that $x^Hx=1.$ Set $r
:=-\Pb(\lambda)x,$ $P_x := I-xx^H$ and $\Lam :=[1, \lam, \ldots,
\lam^m]^T.$
Then  we have \beano \eta_F^{\herm}(\lambda,x,\Pb) &=& \left\{%
\begin{array}{ll}
    \frac{\sqrt{
2\|r\|_{2}^{2}-|x^{H}r|^{2}}}{\|\Lam\|_2} \leq \sqrt{2} \eta(\lam,x,\Pb), & \hbox{if $\lam\in \R,$} \\
    \sqrt{ \|\widehat{r}\|_2^{2}+
\frac{2(\|r\|_{2}^{2}-|x^{H}r|^{2})}{\|\Lam\|_2^2}}, & \hbox{if $\lam \in \C\setminus \R,$} \\
\end{array}%
\right. \\ \eta_2^{\herm}(\lambda,x,\Pb) &=& \left\{%
\begin{array}{ll}
    \eta(\lam,x,\Pb), & \hbox{if $\lam\in \R,$} \\
    \sqrt{ \|\widehat{r}\|_2^{2} +
\frac{\|r\|_{2}^{2}-|x^{H}r|^{2}}{\|\Lam\|_2^2}}, & \hbox{if $\lam \in \C\setminus \R,$} \\
\end{array}%
\right. \eeano where $ \widehat{r} := \bmatrix {\Re (\Lam)^T \\
\Im (\Lam)^T}^\dag \bmatrix{\re (x^Hr) \\ \im (x^Hr)}.$ For the Frobenius norm, define \beano \Delta A_j := \left\{%
\begin{array}{ll}
    \dfrac{\lambda^j}{\|\Lam\|_2^2} (xr^H+rx^H-(r^Hx)xx^H  ) , & \hbox{when $\lam \in \R,$ } \\
    e_j^T \widehat{r} xx^H +
\dfrac{1}{\|\Lam\|_2^2}[\overline{\lam^j}P_xrx^H + \lam^j xr^H
P_x], & \hbox{when $\lam \in \C\setminus \R.$} \\
\end{array}%
\right. \eeano Then  $\Delta \Pb(z):=\sum_{j=0}^m z^j \Delta A_j$
is a unique Hermitian polynomial in $\herm$ such that   $\Pb
(\lambda)x+\Delta \Pb(\lambda)x=0$ and $\nrm {\Delta
\Pb}_F=\eta_F^\herm(\lambda,x,\Pb).$

For the spectral norm, define
$$\Delta A_j :=\left\{%
\begin{array}{ll}
    \dfrac{\lambda^j}{\|\Lam\|_2^2}(rx^H+xr^H-(r^Hx)xx^H ) - \dfrac{\lam^j~x^Hr P_x
rr^H P_x}{\|\Lam\|_2^2~(\|r\|_2^2-|x^Hr|^2)}, & \hbox{when $\lam
\in \R,$}
\\[8pt]  e_j^T
\widehat{r} xx^H + \dfrac{1}{\|\Lam\|_2^2}[\overline{\lam^j}P_x
rx^H + \lam^j xr^H P_x ] - \dfrac{e_j^T \widehat{r} ~P_x
rr^HP_x}{\|r\|_2^2-|x^Hr|^2}, & \hbox{when  $\lam \in \C\setminus \R.$} \\
\end{array}%
\right.$$ Then  $\Delta \Pb(z):=\sum_{j=0}^m z^j \Delta A_j$ is a
Hermitian polynomial in $\herm$ such that   $\Pb (\lambda)x+\Delta
\Pb(\lambda)x=0$ and $\nrm {\Delta
\Pb}_2=\eta_2^\herm(\lambda,x,\Pb).$
\end{theorem}

\noin \pf Again, in view of Theorem~\ref{exist1}, let
 $\Delta \Pb(z)=\sum_{j=0}^m z^j \Delta A_j$ be a Hermitian polynomial such
that $\Delta \Pb(\lam)x+\Pb(\lam)x=0.$ Choosing a unitary matrix
$Q:=[x,~Q_{1}],$ we have
\beano \widetilde{\Delta A_j}:=Q^{H}\Delta A_jQ=\left(%
\begin{array}{cc}
  a_{jj} & a_{j}^{H} \\
  a_{j} &  X_j \\
\end{array}%
\right),~~Q^{H}r=\left(%
\begin{array}{c}
  x^{H}r \\
  Q_{1}^{H}r \\
\end{array}%
\right).\eeano Now
$\Delta \Pb(\lam)x+\Pb(\lam)x=0\Rightarrow \left(%
\begin{array}{c}
  \sum_{j=0}^m \lam^j a_{jj} \\
   \sum_{j=0}^m \lam^j a_{j} \\
\end{array}%
\right)
=\left(%
\begin{array}{c}
  x^{H}r \\
  Q_{1}^{H}r \\
\end{array}%
\right). $ The minimum norm solution of $\sum_{j=0}^m \lam^j
a_{j}=Q_{1}^{H}r$ is given by
$a_j=\frac{\overline{\lam^j}}{\|\Lam\|_2^2}Q_{1}^{H}r.$

Now suppose that $ \lam \in \R.$ Then the minimum norm solution of
$\sum_{j=0}^m \lam^j a_{jj}=x^{H}r$ is given
 by $a_{jj}=\frac{\lam^j}{\|\Lam\|_2^2}x^Hr\in\R.$ Hence for
$\lam\in\R,$ we have
\be\label{bck:hermpoly1} \widetilde{\Delta A_j}=\left(%
\begin{array}{cc}
  \frac{\lam^j}{\|\Lam\|_2^2}x^Hr & (\frac{\lam^j}{\|\Lam\|_2^2}Q_{1}^{H}r)^{H} \\
  \frac{\lam^j}{\|\Lam\|_2^2}Q_{1}^{H}r & X_j \\
\end{array}%
\right), \,\, j = 0:m.\ee
For the Frobenius norm, setting $X_j = 0$ we obtain $\eta_F^\herm
(\lam,x,\Pb)= \frac{\sqrt{
2\|r\|_{2}^{2}-|r^{H}x|^{2}}}{\|\Lam\|_2} $ and the desired Hermitian polynomial
$\Delta \Pb.$

For the spectral norm, setting $ \mu_j :=
\frac{|\lam^j|~\|r\|_2}{\|\Lam\|_2^2}$ and applying
Theorem~\ref{dkw}  to  (\ref{bck:hermpoly1}), we
obtain
$$ X_j = -\frac{\lam^j~x^Hr
(Q_1^Hr)(Q_1^Hr)^H}{\|\Lam\|_2^2~(\|r\|_2^2-|x^Hr|^2)}.$$ This
gives $ \eta^{\herm}(\lambda,x,\Pb) = \frac{\|r\|_{2}}{\|\Lam\|_2}
= \eta(\lam,x,\Pb). $ Now substituting $X_j$ in
(\ref{bck:hermpoly1}) and simplifying the expression,  we obtain
the desired Hermitian polynomial $\Delta \Pb.$

Next, suppose that  $\lambda\in\C\setminus\R.$ Then the minimum
norm solution of $\sum_{j=0}^m \lam^j a_{jj}=x^Hr$ is obtained by
solving
$$\left(%
\begin{array}{c}
  \sum_{j=0}^m \re(\lam^j) a_{jj} \\
  \sum_{j=0}^m \im(\lam^j) a_{jj} \\
\end{array}%
\right)=\left(%
\begin{array}{c}
  \re (x^Hr) \\
  \im (x^Hr) \\
\end{array}%
\right)\Rightarrow \left(%
\begin{array}{c}
  a_{00} \\
  \vdots \\
  a_{mm} \\
\end{array}%
\right)= \left(%
\begin{array}{c}
  \Re (\Lam)^T \\
  \Im (\Lam)^T \\
\end{array}%
\right)^\dag~\left(%
\begin{array}{c}
  \re (x^Hr) \\
  \im (x^Hr) \\
\end{array}%
\right)=: \widehat{r}. $$  Therefore we have $a_{jj}=e_j^T\widehat{r}.$
Hence for $\lam\in \C\setminus\R,$ we have
\be\label{bck:hermpoly2} Q^H \Delta A_jQ=\left(%
\begin{array}{cc}
  e_j^T\widehat{r} & (\frac{\lam^j}{\|\Lam\|_2^2}Q_{1}^{H}r)^{H} \\
  \frac{\lam^j}{\|\Lam\|_2^2}Q_{1}^{H}r & X_j \\
\end{array}%
\right), \,\, j = 0:m.\ee
Thus, for the Frobenius norm, setting $X_j =0$ we obtain
$$\eta_F^\herm(\lam,x,\Pb)= \sqrt{\|\widehat{r}\|_2^{2} + 2
\frac{\|r\|_{2}^{2}-|r^{H}x|^{2}}{\|\Lam\|_2^2}} $$ and
the desired Hermitian polynomial $\Delta \Pb.$

For the spectral norm, setting  $ \mu_j := \sqrt{|e_j^T
\widehat{r} |^2 +
\frac{|\lam^j|^2~(\|r\|_2^2-|x^Hr|^2)}{\|\Lam\|_2^4}}$ and
applying Theorem~\ref{dkw} to the matrix in (\ref{bck:hermpoly2}),
we have
$$ X_j = -\frac{e_j^T \widehat{r}
~(Q_1^Hr)(Q_1^Hr)^H}{\|r\|_2^2-|x^Hr|^2}, \, j=0:m.$$ This gives
$$ \eta_2^{\herm}(\lambda,x,\Pb) = \sqrt{\|\widehat{r}\|_2^{2}
+ \frac{\|r\|^{2}_{2}-|x^Hr|^2}{\|\Lam\|_2^2}}. $$
 Now substituting $X_j$ in
(\ref{bck:hermpoly2}) and simplifying the expression,  we have  $$
\Delta A_j = e_j^T \widehat{r} xx^H +
\frac{1}{\|\Lam\|_2^2}[\overline{\lam^j}P_x rx^H + \lam^j xr^H P_x
] - \frac{e_j^T \widehat{r} ~P_x rr^HP_x}{\|r\|_2^2-|x^Hr|^2}.$$
Hence the results follow. $\blacksquare$

\begin{remark}
If $|x^Hr|=\|r\|_2$ then $\|Q_1^Hr\|_2=0$. Hence considering
$X_j=0, \, j=0:m,$ we obtain the desired results for the spectral
norm.
\end{remark}

Let $ x \in \C^n$ be such that $x^Hx=1.$ If $A \in \C^{n\times n}$
is Hermitian and $A x=0$ then it is easily seen that  $A =
(I-xx^H) Z (I-xx^H)$ for some Hermitian matrix $Z.$ Consequently,
in view of Theorem~\ref{hherm}, we have an analogue of the result
in Corollary~\ref{gpol} for Hermitian matrix polynomials.

Note that, in view of  (\ref{eq:rotb}) and (\ref{eq:rotp}),
structured backward error of $(\lam, x)$ as an approximate
eigenpair of a skew-Hermitian matrix polynomial follows from
Theorem~\ref{hherm}. Indeed, let $\Qb \in \skherm \subset \P$ be a
skew-Hermitian matrix polynomial. Then  $\Pb := i \Qb \in \herm
\subset \P.$ Hence  by (\ref{eq:rotb}) and (\ref{eq:rotp}), we
have $ \eta^\skherm_M(\lam, x, \Qb) = \eta^\herm_M(\lam, x, \Pb).$
Now, let $\Delta \Pb$ be the matrix polynomial given in
Theorem~\ref{hherm} such that $\Pb(\lam)x +\Delta \Pb(\lam) x=0 $
and $\nrm{\Delta \Pb}_M = \eta^\herm_M(\lam, x, \Pb).$  Then
setting $\Delta \Qb := -i \Delta \Pb,$ we  have $\Delta \Qb \in
\skherm$ such that $ \Qb(\lam)x +\Delta\Qb (\lam) x=0$ and $
\nrm{\Delta \Qb}_M = \eta^\skherm_M(\lam, x, \Qb).$

\subsection{H-even and H-odd matrix polynomials}
We now derive structured backward errors of approximate eigenelements of
$H$-even and $H$-odd matrix polynomials. Recall that a matrix
polynomial $\Pb \in \P$ given by $\Pb(z) := \sum^m_{j=0} A_j z^j$
is $H$-even  if and only if $A_j$ is Hermitian  when $j$ is even
(including $j=0$) and $A_j$ is skew-Hermitian when $j$ is odd. Let
$\heven$ and $\hodd,$ respectively, denote the set of $H$-even and
$H$-odd matrix polynomials in $\P.$ Then, as in the case of
Hermitian matrix polynomials in (\ref{eq:rotp}), it is easily seen
that the map \be \label{eq:rotpev} \heven \longrightarrow \hodd,
\,\Pb \longmapsto i \Pb\,\, \mbox{ and } \,\,  \hodd
\longrightarrow \heven, \,\Qb \longmapsto i \Qb \ee are isometric
isomorphisms. Consequently, we only need to prove the results either
for $\heven$ or for $\hodd$ matrix polynomials. Recall that $A^\dag$ is the Moore-Penrose pseudo-inverse of $A$ and $e_j, j=0:m,$ is the canonical basis of $\C^{m+1}.$

\begin{theorem}\label{heven}  Set $\S :=\heven \subset \P.$  Let  $\Pb\in\S$  and $ (\lam, x) \in
\C\times \C^n$ be such that $x^Hx=1.$ Set $r :=-\Pb(\lambda)x,$
$P_x := I-xx^H$ and $\Lam :=[1, \lam, \ldots, \lam^m]^T.$ Then  we
have
\beano \eta_F^{\S} (\lambda,x,\Pb) &=& \left\{%
\begin{array}{ll}
    \frac{\sqrt{
2\|r\|_{2}^{2}-|x^{H}r|^{2}}}{\|\Lam\|_2} \leq \sqrt{2}
\eta(\lam,x,\Pb), & \hbox{if $\lam \in i \,\R,$}
\\[8pt]
    \sqrt{\|\widehat{r}\|_2^{2}+
\frac{2(\|r\|_{2}^{2}-|x^{H}r|^{2})}{\|\Lam\|_2^2}}, & \hbox{if
$\lam\in\C\setminus
i \,\R,$} \\
\end{array}%
\right. \\ \eta_2^{\S} (\lambda,x,\Pb) &=& \left\{%
\begin{array}{ll}
    \eta(\lam,x,\Pb), & \hbox{if $\lam \in i \, \R,$} \\[8pt]
    \sqrt{\| \widehat{r} \|_2^{2}+
\frac{\|r\|^{2}_{2} - |x^Hr|^2}{\|\Lam\|_2^2}}, & \hbox{if $\lam \in \C\setminus i \,\R,$} \\
\end{array}%
\right.\eeano where $\widehat{r} := \bmatrix{\Pi_e \, \Re (\Lam)^T
- (I-\Pi_e)\Im (\Lam)^T \\  \Pi_e \, \Im (\Lam)^T + (I-\Pi_e) \Re
(\Lam)^T}^\dag \bmatrix{\re (x^Hr) \\ \im (x^Hr)}.$ For $j =0:m,$
set $$E_j := \frac{1}{\|\Lam\|_2^2}[\overline{\lam^{j}}P_xrx^H +
\lam^{j} xr^H P_x] \,\,\mbox{ and }\,\, F_j :=
\frac{1}{\|\Lam\|_2^2}[\overline{\lam^{j}}P_x rx^H - \lam^{j} xr^H
P_x ].$$

\noin For the Frobenius norm, define $\dm{\Delta A_j :=
\frac{\overline{\lambda^{j}} }{\|\Lam\|_2^2}
[xr^H+rx^H-(r^Hx)xx^H]}$ when $\lam\in i \,\R$, and
$$\dm{ \Delta A_j :=\left\{%
\begin{array}{ll}
    e_{j}^T \widehat{r} xx^H + E_{j}, & \hbox{if $j$ is even,} \\
    ie_{j}^T \widehat{r} xx^H + F_{j}, & \hbox{if $j$ is odd,} \\
\end{array}%
\right.}$$ when $\lam\in \C\setminus i \,\R,$ for $j =0:m.$ Then
$\Delta \Pb (z) := \sum_{j=0}^m z^j \Delta A_j$ is a unique
$H$-even polynomial in $\S$  such that $\Pb (\lam)x + \Delta \Pb
(\lam)x =0 $ and $\nrm{\Delta \Pb}_F = \eta_F^\S (\lam,x,\Pb).$

For the spectral norm,  define
$$\Delta A_j := \frac{\overline{\lambda^{j}}
}{\|\Lam\|_2^2}[xr^H+rx^H-(r^Hx)xx^H]+\frac{(-1)^{j+1}\lam^{j}~x^Hr
P_x rr^H P_x}{\|\Lam\|_2^2~(\|r\|_2^2-|x^Hr|^2)}$$ when $\lam\in i
\,\R,$ and $$\dm{\Delta A_j :=\left\{%
\begin{array}{ll}
    e_{j}^T \widehat{r} xx^H + E_{j} +
\dfrac{(-1)^{j+1}e_{j}^T \widehat{r} ~P_x rr^H P_x
}{\|r\|_2^2-|x^Hr|^2}, & \hbox{if $j$ is even,} \\ \\
    ie_{j}^T \widehat{r} xx^H + F_{j} -
\dfrac{i \, (-1)^{j+1}  e_{j}^T \widehat{r} ~P_x rr^H P_x
}{\|r\|_2^2-|x^Hr|^2}, & \hbox{if $j$ is odd,} \\
\end{array}%
\right.}$$ when  $\lam\in \C \setminus i \,\R.$  Then $\Delta \Pb
(z) := \sum_{j=0}^m z^j \Delta A_j$ is an $H$-even polynomial in
$\S$ such that $ \Pb (\lam)x + \Delta \Pb (\lam)x =0 $ and
$\nrm{\Delta \Pb}_2 = \eta_2^\S (\lam,x,\Pb).$
\end{theorem}

\noin \pf By Theorem~\ref{exist1} there exists an H-even matrix
polynomial $\Delta \Pb(z)=\sum_{j=0}^m z^j \Delta A_j$ such that
$\Delta \Pb(\lam)=r.$ Now choosing a unitary matrix
$Q:=[x,~Q_{1}],$  we have $\Delta
A_j =  Q\left(%
\begin{array}{cc}
  a_{jj} & a_{j}^{H} \\
  a_{j} & X_{j} \\
\end{array}%
\right)Q^H,$ $X_j^H=X_j $ if $j$ is even, and $\Delta A_j =
    Q\left(%
\begin{array}{cc}
  i a_{jj} & a_{j}^{H} \\
  -a_{j} & Y_{j} \\
\end{array}%
\right)Q^H, \, Y_j^H=-Y_j$ if $j$ is odd. Notice that $a_{jj}$ is
real for all $j.$

Then $\Delta \Pb(\lambda)x=r$ gives
$\left(%
\begin{array}{c}
  \sum_{j\mbox{-even}}
\lam^j a_{jj} + i \sum_{j\mbox{-odd}}  \lam^j a_{jj} \\
  \sum_{j\mbox{-even}} \lam^{j} a_{j} - \sum_{j\mbox{-odd}} \lam^{j} a_{j} \\
\end{array}%
\right)=\left(%
\begin{array}{c}
  x^{H}r \\
  Q^{H}_{1}r \\
\end{array}%
\right).$ The minimum norm solution of $\sum_{j\mbox{-even}}
\lam^{j} a_{j} - \sum_{j\mbox{-odd}} \lam^{j} a_{j}= Q^{H}_{1}r$
is given by $a_{j}=\frac{\overline{\lam^{j}}}{\|\Lam\|_2^2}
Q^{H}_{1}r$ if $j$ is even, and
$a_{j}=-\frac{\overline{\lam^{j}}}{\|\Lam\|_2^2} Q^{H}_{1}r$ if
$j$ is odd.

Now suppose that $\lam\in i \R.$ Then the minimum norm solution
for $a_{jj}$ is given by
$a_{jj}=\frac{\overline{\lam^j}}{\|\Lam\|_2^2} x^Hr$ when $j$ is
even, and $a_{jj}=- \frac{i \, \overline{\lam^j}}{\|\Lam\|_2^2}
x^Hr$ when $j$ is odd. Hence  $a_{jj}\in\R$ when $j$ is even, and
$i a_{jj}\in i \,\R$ when $j$ is odd. Consequently, we have
\be \label{eq:hevenim1} Q^H \Delta A_jQ = \left(%
\begin{array}{cc}
  \frac{\overline{\lam^{j}}}{\|\Lam\|_2^2}
x^Hr & (\frac{\overline{\lam^{j}}}{\|\Lam\|_2^2} Q^{H}_{1}r)^{H}
\\\\
  \frac{\overline{\lam^{j}}}{\|\Lam\|_2^2} Q^{H}_{1}r & X_{j} \\
\end{array}%
\right) \ee when $j$ is even, and \be\label{eq:hevenim2} Q^H \Delta A_j Q = \left(%
\begin{array}{cc}
  \frac{\overline{\lam^{j}}}{\|\Lam\|_2^2}
x^Hr & -(\frac{\overline{\lam^{j}}}{\|\Lam\|_2^2} Q^{H}_{1}r)^{H}
\\\\
  \frac{\overline{\lam^{j}}}{\|\Lam\|_2^2} Q^{H}_{1}r & Y_{j} \\
\end{array}%
\right) \ee when $j$ is odd. Setting $X_{j}=0=Y_{j}$ in
(\ref{eq:hevenim1}) and (\ref{eq:hevenim2}), we obtain $\eta_F^\S
(\lambda,x,\Pb)=\frac{\sqrt{
2\|r\|_{2}^{2}-|r^{H}x|^{2}}}{\|\Lam\|_2}$ and the desired $\Delta
A_j.$

Next, suppose that $\lam\in\C\setminus i \,\R.$  Then
$\sum_{j\mbox{-even}} \lam^j a_{jj} + i \sum_{j\mbox{-odd}} \lam^j
a_{jj}= x^Hr$ gives
$$\left(%
\begin{array}{c}
  \sum_{j\mbox{-even}} \re( \lam^{j}) a_{jj} - \sum_{j\mbox{-odd}} \im( \lam^{j}) a_{jj} \\
  \sum_{j\mbox{-even}} \im( \lam^{j}) a_{jj} + \sum_{j\mbox{-odd}} \re( \lam^{j}) a_{jj} \\
\end{array}%
\right)=\left(%
\begin{array}{c}
  \re (x^Hr) \\
  \im (x^Hr) \\
\end{array}%
\right).$$ Hence we have
\beano \left(%
\begin{array}{c}
  a_{00} \\
  a_{11} \\
  \vdots \\
  a_{mm} \\
\end{array}%
\right) &=&  \left(%
\begin{array}{c}
  \Pi_e \, \Re (\Lam)^T - (I-\Pi_e)\Im (\Lam)^T \\
  \Pi_e \, \Im (\Lam)^T + (I-\Pi_e) \Re (\Lam)^T \\
\end{array}%
\right)^\dag \left(%
\begin{array}{c}
  \re (x^Hr) \\
  \im (x^Hr) \\
\end{array}%
\right) = \widehat{r}   \Rightarrow a_{jj} =  e_j^T
\widehat{r}.\eeano

 Consequently, we have
\be \label{eq:hevenc1} Q^H \Delta A_jQ = \left(%
\begin{array}{cc}
  e_j^T\widehat{r}& (\frac{\overline{\lam^{j}}}{\|\Lam\|_2^2} Q^{H}_{1}r)^{H}
\\\\
  \frac{\overline{\lam^{j}}}{\|\Lam\|_2^2} Q^{H}_{1}r & X_{j} \\
\end{array}%
\right) \ee when $j$ is even, and \be\label{eq:hevenc2} Q^H \Delta A_j Q = \left(%
\begin{array}{cc}
   i e_j^T\widehat{r} & -(\frac{\overline{\lam^{j}}}{\|\Lam\|_2^2} Q^{H}_{1}r)^{H}
\\\\
  \frac{\overline{\lam^{j}}}{\|\Lam\|_2^2} Q^{H}_{1}r & Y_{j} \\
\end{array}%
\right) \ee when $j$ is odd. Now setting $X_{j}=0=Y_{j}$ in
(\ref{eq:hevenc1}) and (\ref{eq:hevenc2}),  we have the desired
matrices $\Delta A_j,$ $j=0:m,$ and $\dm{\eta_F^\S (\lam, x, \Pb)
= \sqrt{\|\widehat{r}\|_2^{2} +2
\frac{\|r\|^{2}_{2}-|x^Hr|^2}{\|\Lam\|_2^2}}}.$ This completes the
proof for the Frobenius norm.

For the spectral norm, consider $ \mu_j :=
\frac{|\lam^{j}|~\|r\|_2}{\|\Lam\|_2^2}$  when $\lambda\in i \R.$
Then applying Theorem~\ref{dkw} to the matrices in
(\ref{eq:hevenim1}) and (\ref{eq:hevenim2}), we obtain
 \beano X_{j} =- \frac{\lam^{j}~x^Hr
(Q_1^Hr)(Q_1^Hr)^H}{\|\Lam\|_2^2~(\|r\|_2^2-|x^Hr|^2)} \mbox{ and
} Y_{j} = \frac{\lam^{j}~x^Hr
(Q_1^Hr)(Q_1^Hr)^H}{\|\Lam\|_2^2~(\|r\|_2^2-|x^Hr|^2)}.\eeano This
gives $ \eta_2^\S(\lam,x,\Pb) = \frac{\|r\|_2}{\|\Lam\|_2}.$ Now
substituting $X_j$ and $Y_j$ in (\ref{eq:hevenim1}) and
(\ref{eq:hevenim2}), we obtain the desired matrices $\Delta A_j, $
$j =0:m.$

When  $\lam\in\C\setminus i \,\R,$ considering  $ \mu_j :=
\sqrt{|e_{j}^T \widehat{r}|^2 +
\frac{|\lam^{j}|^2~(\|r\|_2^2-|x^Hr|^2)}{\|\Lam\|_2^4}}$ and
applying Theorem~\ref{dkw} to the matrices in (\ref{eq:hevenc1})
and (\ref{eq:hevenc2}), we obtain
$$ X_{j} =- \frac{e_{j}^T
\widehat{r}~(Q_1^Hr)(Q_1^Hr)^H}{\|r\|_2^2-|x^Hr|^2} \mbox{ and }
Y_{j} = -\frac{i e_{j}^T
\widehat{r}~(Q_1^Hr)(Q_1^Hr)^H}{\|r\|_2^2-|x^Hr|^2}.
$$ Consequently, we have
$$ \eta_2^{\S}(\lambda,x,\Pb) = \sqrt{\|\widehat{r}\|_2^{2} +
\frac{\|r\|^{2}_{2}-|x^Hr|^2}{\|\Lam\|_2^2}}.$$ Substituting $X_j$
and $Y_j$ in (\ref{eq:hevenc1}) and (\ref{eq:hevenc2}), we obtain
the desired matrices $\Delta A_j, $ $j =0:m.$ $\blacksquare$

Let $ x \in \C$ be such that $x^Hx=1.$ If $ X \in \C^{n\times n}$
is skew-Hermitian and $Xx=0$ then it is easily seen that $X =
(I-xx^H)Z(I-xx^H)$ for some skew-Hermitian matrix $Z.$
Consequently, it follows that an analogue of the result in
Corollary~\ref{cor:teven} holds for $H$-even matrix polynomials.

Observe that, in view of  (\ref{eq:rotb}) and (\ref{eq:rotpev}),
the structured backward error of $(\lam, x)$ as an approximate
eigenpair of an $H$-odd  matrix polynomial follows from
Theorem~\ref{heven}. Indeed, let $\Qb $ be an $H$-odd matrix
polynomial in $\P.$ Set $\S_e := \heven \subset \P$ and $\S_o :=
\hodd \subset \P.$ Then  $\Pb := i \Qb \in \S_e. $ Hence by
(\ref{eq:rotb}) and (\ref{eq:rotpev}), we have $
\eta^{\S_o}_M(\lam, x, \Qb) = \eta^{\S_e}_M(\lam, x, \Pb).$ Now,
let $\Delta \Pb$ be the matrix polynomial given in
Theorem~\ref{heven} such that $ \Delta \Pb \in \S_e,$  $\Pb(\lam)x
+\Delta \Pb(\lam) x=0 $ and $\nrm{\Delta \Pb}_M =
\eta^{\S_e}_M(\lam, x, \Pb).$ Then setting $\Delta \Qb := -i
\Delta \Pb,$ we  have $\Delta \Qb \in \S_o,$  $ \Qb(\lam)x
+\Delta\Qb (\lam) x=0$ and $ \nrm{\Delta \Qb}_M =
\eta^{\S_o}_M(\lam, x, \Pb).$

\subsection{Polynomials with coefficients in Lie and Jordan algebras}
We mention that the structured backward perturbation analysis of
structured matrix polynomials discussed so far can easily be
extended to  more general structured matrix polynomials in which
the coefficient matrices are elements of appropriate Jordan and/or
Lie algebras. Indeed, let $M$ be a unitary matrix such that $M^T=
M$ or $M^T=-M.$ Consider the Jordan algebra $\J := \{A \in
\C^{n\times n} :  M^{-1}A^TM = A\}$ and the Lie algebra $\L :=\{ A
\in \C^{n\times n} : M^{-1}A^TM = -A\}$ associated with the scalar
product $(x,y) \mapsto y^TMx.$ Consider a polynomial $\Pb(z) :=
\sum_{j=0}^m z^j A_j.$ Then by imposing the condition  that  the
polynomial $M\Pb$ given by $M\Pb(z) = \sum_{j=0}^m \lam^j MA_j$ is
either symmetric or skew-symmetric or $T$-even or $T$-odd, we
obtain various structured matrix polynomials. Said differently,
$\S \subset \P$ defines a class of structured matrix polynomials
if $M \S \in \{\sym, \sksym, T\mbox{-}\even, T\mbox{-}\odd\}.$
Hence if $ P \in \S$ then the results obtained in the previous
section are easily extended to  $\Pb $ by replacing $A_j$ and $ r
:=-\Pb(\lam)x$ by $MA_j$ and $ Mr,$ respectively.

Similarly, when $M$ is unitary and $M=M^H$ or $M=-M^H,$ we
consider the Jordan algebra  $\J := \{A \in \C^{n\times n} :
M^{-1}A^HM = A\}$ and the Lie algebra $\L :=\{ A \in \C^{n\times
n} : M^{-1}A^HM = -A\}$ associated with the scalar product $(x,y)
\mapsto y^HMx.$ Then a class of structured matrix polynomials $\S
\subset \P$  is obtained by imposing the condition that  $M \S \in
\{ \herm, \skherm, \heven, \hodd\}.$ Hence the results obtained in
the previous section are easily extended to  $\Pb \in \S$ by
replacing $A_j$ and $ r :=-\Pb(\lam)x$ by $MA_j$ and $ Mr,$
respectively.
In particular, when $M := J,$ where $ J := \left( \ba{cc} 0 & I\\
-I & 0\ea\right) \in \C^{2n\times 2n},$ the Jordan algebra $\J$
consists of skew-Hamiltonian matrices and the Lie algebra $\L$
consists of Hamiltonian matrices. So, for example, considering the
polynomial $\Pb(z) := \sum_{j=0}^m z^j A_j,$ where $A_j$s are
Hamiltonian when $j$ is even and skew-Hamiltonian when $j$ is odd,
we see that the polynomial $J\Pb(z) = \sum_{j=0}^m z^j JA_j$ is
$H$-even. Hence extending the results obtained for $H$-even
polynomial to the case of $\Pb,$ we have the following.

\begin{theorem}\label{hamil} Let $\S$ denote set of polynomials of the form $\Pb(z)=\sum_{j=0}^m z^jA_j$
where $A_{j}$ is Hamiltonian when $j$ is even, and $A_{j}$ is
skew-Hamiltonian when $j$ is odd. Let  $\Pb\in\S$  and $ (\lam, x)
\in \C\times \C^n$ be such that $x^Hx=1.$ Set $r
:=-\Pb(\lambda)x,$ $P_x := I-xx^H$ and $\Lam :=[1, \lam, \ldots,
\lam^m]^T.$ Then we have \beano \eta_F^{\S} (\lambda,x,\Pb) &=& \left\{%
\begin{array}{ll}
    \frac{\sqrt{
2\|r\|_{2}^{2}-|x^{H}Jr|^{2}}}{\|\Lam\|_2} \leq \sqrt{2}
\eta(\lam,x,\Pb), & \hbox{if $\lam \in i \,\R,$}
\\[8pt]
    \sqrt{\|\widehat{r}\|_2^{2}+
\frac{2(\|r\|_{2}^{2}-|x^{H}Jr|^{2})}{\|\Lam\|_2^2}}, & \hbox{if
$\lam\in\C\setminus
i \,\R,$} \\
\end{array}%
\right. \\ \eta_2^{\S} (\lambda,x,\Pb) &=& \left\{%
\begin{array}{ll}
    \eta(\lam,x,\Pb), & \hbox{if $\lam \in i \, \R,$} \\[8pt]
    \sqrt{\| \widehat{r} \|_2^{2}+
\frac{\|r\|^{2}_{2} - |x^HJr|^2}{\|\Lam\|_2^2}}, & \hbox{if $\lam \in \C\setminus i \,\R,$} \\
\end{array}%
\right.\eeano where $\widehat{r} := \bmatrix{\Pi_e \, (\Re
(\Lam)^T) - (I-\Pi_e)(\Im( \Lam)^T) \\  \Pi_e \, (\Im (\Lam)^T) +
(I-\Pi_e) (\Re (\Lam)^T)}^\dag \bmatrix{\re (x^HJr) \\ \im
(x^HJr)}.$
\end{theorem}

\section{Effect of structured linearization on backward error}\label{sec3}
As we have mentioned before, linearization is the standard approach to solving a polynomial
eigenvalue value problem. It is well known that important classes
of structured matrix polynomials admit structured
linearizations~\cite{mackey:thesis,mackey1, mackey2, mackey3}.
However, the process of linearizing a matrix polynomial (structure
preserving or not) has its side effect too. It increases the
sensitivity of eigenvalues of the matrix polynomial~(see,
\cite{higcon, Ba:thesis, BaR09:cond}). Therefore, it is important
to identify linearizations whose eigenelements are almost as
sensitive to perturbations as those of the matrix polynomial.
Obviously, condition numbers of eigenvalues and backward errors of
approximate eigenelements have an important role to play in
identifying such linearizations.  For an unstructured polynomial $
\Pb \in \P,$ Higham et~al.~\cite{higcon,higg1} provide a recipe for choosing a
linearization by analyzing condition numbers of eigenvalues  and
backward errors of approximate eigenpairs. For structured matrix
polynomials,  a recipe for choosing a structured linearization
has been provided in~\cite{BaR09:cond} by analyzing structured condition numbers of eigenvalues. With a view to
identifying optimal and near optimal structured linearizations of a
structured matrix polynomials, in this section, we analyze the
influence of structured linearizations  on the structured backward errors of approximate eigenelements. It turns out that linearizations which minimize structured backward errors also
minimize the structured condition numbers. Therefore our results are consistent with those in~\cite{BaR09:cond}.
We thus provide a recipe for choosing  structured linearizations of a structured
matrix polynomial which minimize structured backward errors as
well as the structured condition numbers.

For a ready reference, we briefly review some basic results about
linearizations of $\Pb \in \P,$ for details, see~\cite{mackey3, mackey:thesis,mackey1, mackey2}.
For our purpose, it is enough to consider the vector space
$\L_1(\Pb)$ given by~\cite{mackey3}
$$ \L_1 (\Pb) := \{ \Lb (\lam) : \Lb (\lam). (\Ls \otimes I_n) = v
\otimes \Pb (\lam), \, v\in \C^m \},$$ where $\Ls :=[\lam^{m-1},
\, \lam^{m-2}, \, \hdots , \, 1]^T$, $\otimes$ is the Kronecker
product and  $v$ is called the right ansatz vector for $\Lb.$
Let $v := [ v_1, \, v_2, \, \hdots, \, v_m]^T\in \C^m$ be a right
ansatz vector. Then the scalar polynomial $\mathsf{p}(x;v) :=
v_1x^{m-1} + v_2x^{m-2} + \hdots + v_{m-1}x + v_m$ is referred to
as the ``$\mathsf{v}$-polynomial" of the vector $v,$
see~\cite{mackey:thesis,mackey3}. The convention is that
$\mathsf{p}(x; v)$ is said to have  a root at $\infty$ whenever $v_1=0.$ Let
$\Lb(\lam) = \lam X + Y \in \L_1 (\Pb)$ be a linearization of
$\Pb$ corresponding to the right ansatz vector $v\in \C^m.$ Then
for $x\in\C^n,$ the following holds
\begin{eqnarray}  \|\Lb(\lam) (\Ls \otimes
x)\|_2 &=& \|v\|_2 \, \|\Pb(\lam)x\|_2, \label{ln1} \\
|(\Ls \otimes x)^T
\Lb(\lam) (\Ls \otimes x)| &=& |\Ls^T v| \, |x^T \Pb(\lam)x|,\label{ln2} \\
|(\Ls \otimes x)^H \Lb(\lam) (\Ls \otimes x)| &=& |\Ls^H v| \,
|x^H \Pb(\lam)x|. \label{ln3}
\end{eqnarray}
Observe from (\ref{ln1}) that $(\lam, x)$ is an eigenelement of $\Pb$
if and only if $(\lam, \, \Ls\otimes x)$ is an eigenelement of
$\Lb.$ Consequently, when $(\lam, x)\in \C\times \C^n$ is
considered as an approximate eigenelement of $\Pb,$ it is natural to
consider $(\lam, \Ls\otimes x) \in \C\times \C^{mn}$ as an
approximate eigenelement of $\Lb$ and vice-versa. We denote the (unstructured)
backward error of $(\lam, \Ls\otimes x)$ as an approximate eigenelement of $\Lb$ by $\eta(\lam,\Ls\otimes
x,\Lb;v)$ so as to show the dependence of the backward error on
the ansatz vector $v.$ Similarly, we denote the structured
backward by $\eta^{\S}_M(\lam,\Ls\otimes
x,\Lb;v),$ where $ M \in \{2, F\}.$

Now suppose that $(\lam, x) \in \C\times \C^n$ with $x^Hx=1$ is an approximate
eigenpair of $\Pb.$ In view of (\ref{ln1}) - (\ref{ln3}), we only need to
consider ansatz vectors $v$ having unit norm. We use the inequality \be\label{bckerr:insig} \sqrt{\frac{m+1}{2m}} \leq \frac{
\|\Lam\|_2 }{\|\Ls\|_2 \, \|(\lam, \, 1)\|_2} \leq 1 \ee which is derived in (Lemma A.1, \cite{higcon}).

\begin{theorem}\label{unstr:bck}
Let $\Pb\in\P$ be  regular and $\Lb\in\L_1(\Pb)$ be a linearization of $\Pb$ corresponding to the normalized
right ansatz vector $v.$  Let $(\lam,x) \in \C\times
\C^n$ be such that $x^Hx =1.$ Set $\Ls :=[ \lam^{m-1}, \ldots, \lam,
1]^T.$ Then we have $$\sqrt{\frac{m+1}{2m}} \leq \frac{\eta (\lam,\Ls
\otimes x,\Lb;v)}{\eta (\lam,x,\Pb)} \leq 1.$$
\end{theorem}

\noin\pf By~(\ref{unstr:bckerr}) and (\ref{ln1}), we have \beano
\eta (\lam,\Ls \otimes x,\Lb;v) &=& \frac{\|\Lb(\lam)(\Ls \otimes
x)\|_2}{\|(\Ls \otimes x)\|_2 \, \|(\lam, \, 1)\|_2} =
\frac{\|v\|_2
\|\Pb(\lam)x\|_2}{\|(\Ls \otimes x)\|_2 \, \|(\lam, \, 1)\|_2}\\
&=& \frac{\|v\|_2 \|\Lam\|_2 \|x\|_2}{\|(\Ls \otimes x)\|_2 \,
\|(\lam, \, 1)\|_2} \eta (\lam,x,\Pb) \\ &=& \frac{ \|\Lam\|_2
}{\|\Ls\|_2 \, \|(\lam, \, 1)\|_2} \eta (\lam,x,\Pb). \eeano Hence
 by~(\ref{bckerr:insig}) the desired result follows. $\blacksquare$

Theorem~\ref{unstr:bck} shows that as far as the backward
errors of approximate eigenelements of $\Pb$ are concerned, any
linearization from $\L_1(\Pb)$ is as good as any other provided that the linearization is associated to a normalized right ansatz vector. In contrast, restricting $\Lb$ in $\D\L(\Pb)$ (see, \cite{mackey3}), it is shown in \cite{BaR09:cond} that the condition number of an eigenvalues $\lam$ of $\Pb$ is increased at least by $\delta(\lam, v)$ and at most by $\sqrt{2} \, \delta(\lam, v),$ where $\delta(\lam, v) := \|\Ls\|_2/|\mathsf{p}(x;v)|,$ see also~\cite{higcon}.

For a structured matrix polynomials, there exists infinitely many structured linearizations, see~\cite{mackey1, mackey2, mackey:thesis}. For the structures we consider in this paper, we consider structured linearization from  $\L_1(\Pb).$ For a ready reference, we summarize in Table~\ref{tablebckerr:ansatzvectors} the condition on ansatz vector  for a structured linearization,  see~\cite{mackey1,mackey2}. The matrix $\Sigma$ in Table~\ref{tablebckerr:ansatzvectors} is given by $ \Sigma = \text{diag}\{
  (-1)^{m-1},  (-1)^{m-2},\ldots,(-1)^0\}.$

\begin{table}[h]
\begin{center}\renewcommand{\arraystretch}{1.2}
\begin{tabular}{|r|c|l|}
  \hline
  $\S$ & Structured Linearization & ansatz vector \\
  \hline
  \hline
  $\sym $ & sksymm & $v \in \C^m$ \\
  \hline
   $\sksym$ & skew-symm & $v \in \C^m$ \\
   \hline
  $T$-$\even$ & $T$-even & $\Sigma v = v$ \\ \cline{2-3}
   & $T$-odd & $\Sigma v =-v$ \\
  \hline
  $T$-$\odd$ & $T$-even & $\Sigma v =-v$ \\ \cline{2-3}
  & $T$-$\odd$ & $\Sigma v = v$ \\
    \hline
  \hline
 $\herm$ & Herm & $v \in \R^m $ \\ \cline{2-3}
  & skew-Herm & $v \in i\R^m $ \\
  \hline
$\skherm$ & Herm & $v \in i\R^m $ \\ \cline{2-3}
 & skew-Herm & $v \in \R^m $ \\
  \hline
  $\heven$ & $H$-even & $\Sigma v = \overline{v}$ \\ \cline{2-3}
   & $T$-odd & $\Sigma v = -\overline{v}$ \\
  \hline
  $\hodd$ & $H$-even & $\Sigma v = -\overline{v}$ \\ \cline{2-3}
  & $H$-odd & $\Sigma v = \overline{v}$ \\
  \hline
\end{tabular}
\caption{\label{tablebckerr:ansatzvectors} Admissible ansatz vectors for structured linearizations.}
\end{center}
\end{table}

Recall that $\eta(\lam,x,\Pb) \leq \eta_M^\S(\lam,x, \Pb).$  Similarly, for a
structured linearization from $\L_1(\Pb)$ we have $\eta(\lam,\Ls\otimes
x,\Lb;v) \leq \eta_M^\S(\lam,\Ls\otimes x, \Lb;v),$ where $v$ is the
ansatz vector. With a view to understanding the effect of structure preserving linearizations on the backward errors of approximate eigenelements of structured matrix polynomials,  in this section we compare $\eta(\lam, x, \Pb)$ and $ \eta^\S_M(\lam, x, \Pb)$ with $ \eta^\S_M(\lam, x, \Lb).$

\begin{corollary}\label{unbck:corollary}
Let $\Pb\in\S$ and $\Lb\in\L_1(\Pb)$ be a structured
linearization  corresponding to the normalized ansatz vector $v.$ Then for $M\in \{2, F\},$ we have
$$\frac{\eta^\S_M(\lam,\Ls\otimes x,\Lb;v)}{\eta(\lam,x,\Pb)} \geq
\sqrt{\frac{m+1}{2m}}.$$
\end{corollary}

\noin\pf By Theorem~\ref{unstr:bck} we have
$$\frac{\eta^\S_M(\lam,\Ls\otimes x,\Lb;v)}{\eta(\lam,x,\Pb)} \geq \frac{\eta(\lam,\Ls\otimes x,\Lb;v)}{\eta(\lam,x,\Pb)} \geq
\sqrt{\frac{m+1}{2m}}.$$ Hence the proof. $\blacksquare$

\subsection{Symmetric and skew-symmetric linearizations}
For a symmetric matrix polynomial $\Pb$, any ansatz vector $v$
yields a potential symmetric linearization. Recall that an ansatz vector $v$ is always assumed to be normalized, that is, $\|v\|_2 =1.$  We now show that structure preserving linearizations of symmetric and skew-symmetric  matrix polynomials have almost no adverse effect on the backward errors of approximate eigenelements.

\begin{theorem}\label{bck:compsym1}
Let $\S$ be the space of symmetric matrix polynomials and $\Pb \in \S.$
Let $\Lb \in \L_1(\Pb)$ be a symmetric linearization of $\Pb$ with normalized ansatz
vector $v.$ Finally, let $(\lam,x)\in\C \times \C^n$ be such that
$\|x\|_2=1.$ Then we have
\beano \sqrt{\frac{m+1}{2m}} &\leq &
\dfrac{ \eta^\S_F (\lam, \Ls\otimes x,\Lb; v)}{ \eta^\S_F(\lam, x, \Pb)} \leq \dfrac{
\eta^\S_F (\lam, \Ls\otimes x,\Lb; v)}{ \eta(\lam, x, \Pb)} \leq \sqrt{2},\\ \sqrt{\frac{m+1}{2m}} &\leq&
\dfrac{\eta^\S_2 (\lam, \Ls\otimes x,\Lb; v)}{\eta^\S_2(\lam, x, \Pb)}= \dfrac{\eta (\lam, \Ls\otimes x,\Lb; v)}{\eta(\lam, x, \Pb)}  \leq 1. \eeano
\end{theorem}

\noin\pf For the Frobenius norm, by Theorem~\ref{tsym} we have
$$\dfrac{\eta^\S_F(\lam,\Ls\otimes x,\Lb;v)}{\eta^\S_F(\lam,x,\Pb)}
= \dfrac{\sqrt{2\|r\|_2^2-\frac{|\Ls^Tv|^2}{\|\Ls\|_2^2}|x^Tr|^2}}{\sqrt{2\|r\|_2^2
-|x^Tr|^2 }} \cdot \dfrac{\|\Lam\|_2}{\|\Ls\|_2\|(\lam, \,
1)\|_2}, $$ where  $ r:=-\Pb(\lam)x.$ Hence by
(\ref{bckerr:insig}) we have $\dfrac{\eta^\S_F(\lam,\Ls\otimes
x,\Lb;v)}{\eta^\S_F(\lam,x,\Pb)}\geq \sqrt{\frac{m+1}{2m}}.$
Next, since  $\|r\|_2 \leq \sqrt{2
\|r\|_2^2 - \frac{|\Ls^Tv|^2}{\|\Ls\|_2^2} |x^Tr|^2} \leq \sqrt{2} \,
\|r\|_2,$ we have

$$\dfrac{\eta^\S_F(\lam,\Ls\otimes
x,\Lb;v)}{\eta^\S_F(\lam,x,\Pb)} \leq
\dfrac{\eta^\S_F(\lam,\Ls\otimes x,\Lb;v)}{\eta(\lam,x,\Pb)} \leq \sqrt{2} \dfrac{\|\Lam\|_2}{\|\Ls\|_2\|(\lam, \,
1)\|_2} \leq \sqrt{2}.$$

Finally, by Theorem~\ref{tsym}, we have structured  and unstructured backward errors are the same for the spectral norm. Hence the desired results follow from Theorem~\ref{unstr:bck}. $\blacksquare$

For skew-symmetric linearizations of skew-symmetric matrix polynomials, we have the following result.

\begin{theorem}\label{bck:compsksym1}
Let $\S$ be the space of skew-symmetric matrix polynomials and $\Pb \in \S.$
Let $\Lb \in \L_1(\Pb)$ be a skew-symmetric linearization of $\Pb$ with normalized
ansatz vector $v.$ Finally, let $(\lam,x)\in\C \times \C^n$ be such that
$\|x\|_2=1.$ Then for $ M \in \{2, F\}$ we have
$$ \sqrt{\frac{m+1}{2m}} \leq
\dfrac{\eta^\S_M (\lam, \Ls\otimes x,\Lb; v)}{\eta^\S_M(\lam, x, \Pb)}= \dfrac{\eta(\lam, \Ls\otimes x,\Lb; v)}{\eta(\lam, x, \Pb)}  \leq 1.$$
\end{theorem}

\noin\pf  By Theorem~\ref{tskewsym} we
have $$\dfrac{\eta^\S_M (\lam, \Ls\otimes x,\Lb; v)}{\eta^\S_M(\lam, x, \Pb)}= \dfrac{\eta(\lam, \Ls\otimes x,\Lb; v)}{\eta(\lam, x, \Pb)}.$$ Hence desired result follows from
Theorem~\ref{unstr:bck}.
$\blacksquare$

Thus we conclude that for a symmetric/skew-symmetric matrix polynomial a structure preserving linearization automatically ensures that the  backward errors of approximate eigenelements are least affected by the conversion the polynomial eigenvalue problem into a generalized eigenvalue problem of larger dimension. Moreover, as shown in~\cite{BaR09:cond} this choice also ensures that the linearization has a mild influence on the structured condition numbers of eigenvalues of the polynomial.

\subsection{$T$-even and $T$-odd linearizations}

Now we analyze  $T$-even and $T$-odd linearizations.  Note that a $T$-even (resp., $T$-odd) polynomial admits $T$-even as well as $T$-odd linearizations which preserve the spectral symmetry of the  $T$-even (resp., $T$-odd) polynomial.

\begin{theorem}\label{bckerr:compteven} Let $\Pb\in \P$ be a $T$-even polynomial and $ (\lam, x) \in \C\times \C^n$ be such that $ \|x\|_2 =1.$  Let $\S_e \subset \L_1(\Pb)$ and $\S_o \subset \L_1(\Pb),$ respectively,  denote the space  of $T$-even and $T$-odd pencils. Finally, let $\Lb_e \in \S_e$ (resp.,
$\Lb_o \in \S_o$) be $T$-even (resp. $T$-odd) linearization of $\Pb$ with normalized ansatz vector  $v =\Sigma
v$ (resp., $v=- \Sigma v$). Then  for $M \in \{2, F\}$ we have the following.
\begin{enumerate} \item If $|\lam|\leq 1$ then
$\sqrt{\frac{m+1}{2m}} \leq
\dfrac{\eta_M^{\S_e} (\lam,\Ls \otimes x,\Lb_e; v)}{\eta(\lam,x,\Pb)}
\leq \sqrt{2}.$

\item  If $|\lam|\geq 1$ then $\sqrt{\frac{m+1}{2m}} \leq
\dfrac{\eta_M^{\S_o} (\lam,\Ls \otimes x,\Lb_o; v)}{\eta(\lam,x,\Pb)}
\leq \sqrt{2}.$\end{enumerate}
\end{theorem}

\noin\pf First consider the $T$-even linearization $\Lb_e.$ Then by Theorem~\ref{teven} we have  \be\label{fro:evnl} \frac{\eta^{\S_e}_F
(\lam,\Ls \otimes x,\Lb_e; v)}{\eta(\lam, x, \Pb)}
= \frac{\left(\sqrt{2 \, \|r\|_2^2 + \frac{(|\lam|^2
-1)|\Ls^T v|^2}{\|\Ls\|_2^2} \, |x^Tr|^2}\right) \|\Lam\|_2}{\|r\|_2 \, \|\Ls\|_2  \,
\|(\lam, \, 1)\|_2}, \ee where $r := -\Pb(\lam)x.$

Now for $|\lam|\leq 1,$  we have $ \|r\|_2 \leq \sqrt{2
\|r\|_2^2 + \frac{(|\lam|^2-1) |\Ls^Tv|^2}{\|\Ls\|_2^2} |x^Tr|^2} \leq
\sqrt{2}\, \|r\|_2.$ Hence by~(\ref{bckerr:insig}) we obtain the desired results for
the Frobenius norm.

Again by Theorem~\ref{teven}, we have
\be\label{2:evnl} \frac{ \eta^{\S_e}_2 (\lam, \Ls \otimes x,\Lb_e; v)}{\eta(\lam, x, \Pb)}
 = \frac{\left(\sqrt{ \|r\|_2^2 + \frac{|\lam|^2 \, |\Ls^T
v|^2}{\|\Ls\|_2^2} \, |x^Tr|^2}\right)\|\Lam\|_2}{ \|r\|_2 \|\Ls\|_2 \, \|(\lam ,
\, 1)\|_2}. \ee Notice that  $\|r\|_2 \leq \sqrt{
\|r\|_2^2 +|\lam|^2\frac{|\Ls^Tv|^2}{\|\Ls\|_2^2} |x^Tr|^2} \leq
\sqrt{1+|\lam|^2}\, \|r\|_2$  for $\lam\in\C.$ Hence by~(\ref{bckerr:insig}) we obtain the desired
result for the spectral norm.

Next, consider the $T$-odd linearization $\Lb_o.$ Then by Theorem~\ref{todd}  we have
\be\label{fro:odl}
\frac{\eta^{\S_o}_F (\lambda,\Ls \otimes x, \Lb_o;v)}{\eta(\lam, x, \Pb)} =
\frac{\left(\sqrt{
2 \|r\|_2^2 + (\frac{1}{|\lam|^2}-1)\frac{|\Ls^T
v|^2}{\|\Ls\|_2^2} |x^Tr|^2}\right) \|\Lam\|_2}{ \|r\|_2 \|\Ls\|_2  \,
\|(1, \lam)\|_2}, \ee for $\lam \neq 0.$ Now for $|\lam|\geq 1,$ we have  $$  \|r\|_2
\leq \sqrt{ 2 \|r\|_2^2 +(|\lam|^{-2}-1)\frac{|\Ls^Tv|^2}{\|\Ls\|_2^2}
|x^Tr|^2} \leq \sqrt{ 2} \|r\|_2.$$ Hence by (\ref{bckerr:insig})we obtain the desired result for the Frobenius norm.

Again by Theorem~\ref{todd} we have \be\label{2:odl}
\frac{\eta^{\S_o}_2 (\lambda,\Ls \otimes x, \Lb_o;v)}{\eta(\lam, x, \Pb)}=
\frac{\left(\sqrt{\|r\|_2^2  + \frac{|\Ls^T v|^2}{|\lam|^2 \, \|\Ls\|_2^2}
|x^Tr|^2}\right)\|\Lam\|_2}{ \|r\|_2\, \|\Ls\|_2  \|(1, \lam)\|_2} \ee for $\lam \neq 0.$
For  $|\lam| \geq 1,$ we have $ \|r\|_2
\leq \sqrt{ \|r\|_2^2 + |\lam|^{-2}\frac{|\Ls^Tv|^2}{\|\Ls\|_2^2}
|x^Tr|^2} \leq \sqrt{2}
\|r\|_2.$  Hence by~(\ref{bckerr:insig}) we obtain the desired
result follows for the spectral norm.$\blacksquare$

\begin{remark}
We mention that the bounds in Theorem~\ref{bckerr:compteven} also hold when $\Pb$ is $T$-odd  with the role of $T$-even and $T$-odd linearizations are reversed, that is,  by interchanging the role of $\Lb_e$ and $\Lb_o$ we obtain the desired bounds.
\end{remark}

Next, comparing  $\eta^\S_2(\lam, x, \Pb)$ with $\eta^\S_2(\lam, x, \Lb, v)$ we have the following result.

\begin{theorem} Suppose that the assumptions of Theorem~\ref{bckerr:compteven} hold. Let $\S \subset\P$ denote the set of $T$-even polynomials. Then we have \begin{enumerate} \item If $|\lam|\leq 1:$
$\sqrt{\frac{m+1}{4m}} \leq\dfrac{\eta^{\S_e}_2(\lam,\Ls\otimes
x,\Lb_e;v)}{\eta^\S_2(\lam,x,\Pb)}\leq \sqrt{2}.$ \item If $|\lam|\geq 1:$  $ \sqrt{\frac{m+1}{4m}} \leq\dfrac{\eta^{\S_o}_2(\lam,\Ls\otimes
x,\Lb_o;v)}{\eta^\S_2(\lam,x,\Pb)}\leq \sqrt{2},$ when $m$ is even and\\
$\sqrt{\frac{m+1}{2m}} \frac{1}{\|(1, \,\, \lam)\|_2}\leq\dfrac{\eta^{\S_o}_2(\lam,\Ls\otimes
x,\Lb_o;v)}{\eta^\S_2(\lam,x,\Pb)}\leq \sqrt{2},$ when $m$ is odd.
\end{enumerate}
\end{theorem}

\noin\pf Note that the upper bounds follow from Theorem~\ref{bckerr:compteven}. We now derive the lower bounds.

First suppose that $|\lam|\leq 1.$  Then it is easy to see that $\|(I-\Pi_e)(\Lam)\|_2 \leq \|\Pi_e(\Lam)\|_2.$ Hence by Theorem \ref{teven} we have $$ \eta^\S_2 (\lam,x,\Pb) \leq \frac{\|r\|_2}{\|\Lam\|_2} \sqrt{1 + \frac{\|(I-\Pi_e)(\Lam)\|_2^2}
{\|\Pi_e(\Lam)\|_2^2}} \leq \frac{\sqrt{2}\|r\|_2}{\|\Lam\|_2}. $$ On the other hand, by (\ref{2:evnl}) we have
$\dm{\eta_2^{\S_e}(\lam, \Lambda_{m-1} \otimes x, \Lb_e; v ) \geq \frac{\|r\|_2}{\|\Lambda_{m-1}\|_2 \|(1,\lam)\|_2}.}$
Consequently, by (\ref{bckerr:insig}) we have $$\frac{\eta_2^{\S_e}(\lam, \Lambda_{m-1} \otimes x, \Lb_e;v)}{\eta^\S_2 (\lam,x,\Pb)} \geq
\frac{\|\Lam\|_2}{\sqrt{2}\|\Lambda_{m-1}\|_2 \|(1.\lam)\|_2} \geq \frac{1}{2}\sqrt{\frac{m+1}{m}}.$$

Next suppose that $|\lam|\geq1$ and consider the $T$-odd linearization $\Lb_o.$  Then it is easy to check that
$\|(I-\Pi_e)(\Lam)\|_2 \leq \|\Pi_e(\Lam)\|_2$ when $m$ is even and the desired result follows by similar arguments as above. Now suppose that $m$ is odd. Then it is easy to see that $\|(I-\Pi_e)(\Lam)\|_2^2 = |\lam|^2 \|\Pi_e(\Lam)\|_2^2.$ Hence  by Theorem \ref{teven} we have $$ \eta^\S_2 (\lam,x,\Pb) \leq \frac{\|r\|_2}{\|\Lam\|_2} \sqrt{1 + \frac{\|(I-\Pi_e)(\Lam)\|_2^2}
{\|\Pi_e(\Lam)\|_2^2}} \leq \frac{\sqrt{1 +|\lam|^2}\|r\|_2}{\|\Lam\|_2}. $$ Further by (\ref{2:odl}) we have
$\dm{\eta_2^{\S_o}(\lam, \Lambda_{m-1} \otimes x, \Lb_o; v ) \geq \frac{\|r\|_2}{\|\Lambda_{m-1}\|_2 \|(1,\lam)\|_2}.}$ Hence by (\ref{bckerr:insig}) we have
$$\frac{\eta_2^{\S_o}(\lam, \Lambda_{m-1} \otimes x, \Lb_o;v)}{\eta^\S_2 (\lam,x,\Pb)} \geq
\frac{\|\Lam\|_2}{\|\Lambda_{m-1}\|_2 \|(1.\lam)\|_2^2} \geq \frac{1}{\sqrt{2}\|(1, \,\, \lam )\|_2}\sqrt{\frac{m+1}{m}}.$$  This completes the proof. $\blacksquare$

For $T$-odd polynomials, we have the following result.

\begin{theorem}
Let $\S\subset \P$ denote the space of $T$-odd matrix polynomials and $ \Pb \in \S.$ Let $\S_e \subset \L_1(\Pb)$ and $\S_o \subset \L_1(\Pb),$ respectively,  denote the space  of $T$-even and $T$-odd pencils. Finally, let $\Lb_e \in \S_e$ (resp., $\Lb_o \in \S_o$) be $T$-even (resp. $T$-odd) linearization of $\Pb$ with normalized ansatz vector  $v =\Sigma v$ (resp., $v=- \Sigma v$).  Then for $(\lam,x)\in\C \times \C^n$ with $\|x\|_2 =1,$ we have the following.

\begin{enumerate} \item If $|\lam|\leq 1:$ $\sqrt{\frac{m+1}{6m}} \leq\dfrac{\eta^{\S_o}_2(\lam,\Ls\otimes
x,\Lb_o;v)}{\eta^\S_2(\lam,x,\Pb)}\leq 1.$ \item If $|\lam|\geq 1:$  $\dfrac{1}{\|(\sqrt{2}, \,\, \lam)\|_2} \sqrt{\frac{m+1}{m}} \leq\dfrac{\eta^{\S_e}_2(\lam,\Ls\otimes
x,\Lb_e;v)}{\eta^\S_2(\lam,x,\Pb)}\leq 1,$ when $m$ is even and\\
$\sqrt{\frac{m+1}{4m}}\leq\dfrac{\eta^{\S_o}_2(\lam,\Ls\otimes
x,\Lb_o;v)}{\eta^\S_2(\lam,x,\Pb)}\leq 1,$ when $m$ is odd.
\end{enumerate}
\end{theorem}

\noin\pf By Theorem~\ref{todd} we have $\dm{\eta_2^\S(\lam,x;\Pb) = \frac{1}{\|\Lam\|_2}\sqrt{\|r\|_2^2 + \frac{\|\Pi_e(\Lam)\|_2^2}{\|(I-\Pi_e)(\Lam)\|_2^2}|x^Tr|^2}.}$ It is easy to see that \be\label{rel:1} \frac{|\lam|^2}{1+|\lam|^2} \geq \frac{\|(I-\Pi_e)(\Lam)\|_2^2}{\|\Lam\|_2^2}\ee with equality holds for odd $m.$
Now, by (\ref{2:odl}) we have $$\eta_2^{\S_o}(\lam,\Ls\otimes x,\Lb_o;v) = \frac{1}{\|\Ls\|_2 \|(1, \,\, \lam)\|_2} \sqrt{\|r\|_2^2 + |\lam|^{-2} \frac{|\Ls^Tv|^2}{\|\Ls\|_2^2} |x^Tr|^2}.$$ Since by (\ref{rel:1}), $\frac{\|\Pi_e(\Lam)\|_2^2}{\|(I-\Pi_e)(\Lam)\|_2^2}\geq |\lam|^{-2}$ with equality holds for odd $m,$   we have \be\label{re:odd}\frac{\eta^{\S_o}_2(\lam,\Ls\otimes
x,\Lb_o;v)}{\eta^\S_2(\lam,x;\Pb)}\leq 1.\ee  Also it is easy to check that $\|(I-\Pi_e)(\Lam)\|_2 \leq \|\Lam\|_2 \leq \sqrt{2} \|\Pi_e(\Lam)\|_2 $ whenever $|\lam|\leq 1.$ Consequently we have $\frac{\|(I-\Pi_e)(\Lam)\|_2}{\|\Pi_e(\Lam)\|_2} \leq \sqrt{2}.$ This yields $$\frac{\eta^{\S_o}_2(\lam,\Ls\otimes
x,\Lb_o;v)}{\eta^\S_2(\lam,x;\Pb)} \geq  \sqrt{\frac{m+1}{6m}}.$$

 Next suppose that $|\lam| \geq 1.$ If $m$ is even then its obvious that $\frac{\|\Pi_e(\Lam)\|_2^2}{\|(I-\Pi_e)(\Lam)\|_2^2} \geq |\lam|^2.$ Hence by (\ref{2:evnl}) and (\ref{bckerr:insig}) we have $$\frac{\eta^{\S_e}_2(\lam,\Ls\otimes
x,\Lb_e;v)}{\eta^\S_2(\lam,x;\Pb)} \leq \frac{\|\Lam\|_2}{\|\Ls\|_2\|(1, \, \, \lam)\|_2} \frac{\sqrt{\|r\|_2^2 + |\lam|^2\frac{|\Ls^Tv|^2}{\|\Ls\|_2^2}|x^Tr|^2}}{\sqrt{\|r\|_2^2 + |\lam|^2|x^Tr|^2}} \leq 1. $$ Further  using the fact $\dfrac{\|\Pi_e(\Lam)\|_2^2}{\|(I-\Pi_e)(\Lam)\|_2^2} \leq 1+|\lam|^2$ we have $$\frac{\eta^{\S_e}_2(\lam,\Ls\otimes x,\Lb_e;v)}{\eta^\S_2(\lam,x;\Pb)} \geq \frac{1}{\|(\sqrt{2}, \,\, \lam)\|_2}\sqrt{\frac{m+1}{m}}.$$ On the other hand, if $m$ is odd and $|\lam|\geq 1$  then the lower bound follows from the fact that $\dfrac{\|\Pi_e(\Lam)\|_2^2}{\|(I-\Pi_e)(\Lam)\|_2^2} = \frac{1}{|\lam|^2} \leq 1.$ This completes the proof. $\blacksquare$

The moral of the story is that for computing eigenelements of a $T$-even matrix polynomial $\Pb,$  it is advisable to solve $T$-even as well as $T$-odd linearizations of $\Pb$ and then choose a computed eigenpair $(\lam, x)$  from $T$-even or $T$-odd linearization according as $|\lam| \leq 1$ or $|\lam| \geq 1.$ In contrast, when $P$ is $T$-odd it is advisable to choose  $(\lam, x)$  from $T$-even linearization only when $|\lam|\geq 1$ and the degree of $\Pb$ is even, otherwise choose $(\lam, x)$ from $T$-odd linearization of $\Pb.$ This choice ensures that the linearizations have almost no adverse effect  on the backward error of the computed eigenelement $(\lam, x).$ We arrived at the same conclusion in ~\cite{BaR09:cond} by analyzing the effect of structure preserving linearizations on the structured condition numbers of eigenvalues of the polynomial $\Pb.$

\subsection{Hermitian and $H$-even linearizations}

First, we consider Hermitian matrix polynomials.  Note that a
Hermitian matrix polynomial admits  Hermitian and
skew-Hermitian linearizations both preserving the spectral symmetry of
the Hermitian polynomial.

\begin{theorem}\label{hherm2}
Let $\Pb \in \P$ be Hermitian and $(\lam, x) \in \C\times \C^n$ be such that $\|x\|_2 =1.$ Let $ \S \in \{\herm, \skherm\}$ and $ \Lb \in \S $ be a linearization of $\Pb$ with normalized ansatz vector $v.$ If $\lam \in \R$ then  we have $$\sqrt{\frac{m+1}{2m}} \leq
\dfrac{\eta^{\S}_F(\lam,\Ls\otimes x,\Lb; v)}{\eta^{\herm}_F(\lam,x,\Pb)} \leq
\sqrt{2}  \mbox{ and } \sqrt{\frac{m+1}{2m}} \leq
\dfrac{\eta^{\S}_2(\lam,\Ls\otimes x,\Lb;
v)}{\eta^{\herm}_2(\lam,x,\Pb)} \leq 1.$$
The same bounds hold when $\Pb$ is skew-Hermitian.
\end{theorem}

\noin\pf First, suppose that $ \S = \herm$ so that $\Lb$ is a Hermitian linearization of $\Pb.$ For $\lam \in \R,$
by Theorem~\ref{hherm}, we have
$$\dfrac{\eta^\S_F(\lam,\Ls\otimes x,\Lb;v)}{\eta^{\herm}_F(\lam,x,\Pb)}
= \dfrac{\sqrt{2\|r\|_2^2-\frac{|\Ls^Hv|^2}{\|\Ls\|_2^2}|x^Hr|^2}}{\sqrt{2\|r\|_2^2
-|x^Hr|^2 }} \cdot \dfrac{\|\Lam\|_2}{\|\Ls\|_2\|(\lam, \,
1)\|_2}, $$ where  $ r:=-\Pb(\lam)x.$ Hence by
(\ref{bckerr:insig}) we have $\dfrac{\eta^\S_F(\lam,\Ls\otimes
x,\Lb;v)}{\eta^{\herm}_F(\lam,x,\Pb)}\geq \sqrt{\frac{m+1}{2m}}.$ Next,
since  $\|r\|_2 \leq \sqrt{2
\|r\|_2^2 - \frac{|\Ls^Hv|^2}{\|\Ls\|_2^2} |x^Hr|^2} \leq \sqrt{2} \,
\|r\|_2,$ we have

$$\dfrac{\eta^\S_F(\lam,\Ls\otimes
x,\Lb;v)}{\eta^{\herm}_F(\lam,x,\Pb)} \leq
\dfrac{\eta^\S_F(\lam,\Ls\otimes x,\Lb;v)}{\eta(\lam,x,\Pb)} \leq \sqrt{2} \dfrac{\|\Lam\|_2}{\|\Ls\|_2\|(\lam, \,
1)\|_2} \leq \sqrt{2}.$$

For the spectral norm, by Theorem~\ref{hherm}, structured  and unstructured backward errors are the same when $\lam \in \R.$ Hence the desired results follow from Theorem~\ref{unstr:bck}.

Finally, since the backward errors are the same for Hermitian and skew-Hermitian pencils, the above hounds obviously hold for the case when $ \S= \skherm.$  $\blacksquare$

This shows that a structured linearization of a Hermitian matrix polynomial does not have adverse effect on the backward errors of approximate eigenelements when the approximate eigenvalues are real. On the other hand, when the approximate eigenvalues are complex, the structured backward errors are not amenable to easy comparisons.
Indeed, under the assumptions of Theorem~\ref{hherm2}, when $\lam \in \C\setminus\R$ by Theorem~\ref{hherm} a little calculation shows that
$$ \dfrac{\eta^\S_F(\lam,\Ls\otimes
x,\Lb; v)}{\eta(\lam, x, \Pb)} \leq \sqrt{2 + \frac{\|\widehat{r}\|_2^2}{\|r\|_2^2}} \mbox{ and }
\dfrac{\eta^\S_2(\lam,\Ls\otimes
x,\Lb; v)}{\eta(\lam, x, \Pb)} \leq \sqrt{1 + \frac{\|\widehat{r}\|_2^2}{\|r\|_2^2}},$$
where $ r := -\Pb(\lam)x$ and $ \widehat{r} = r_h:= \bmatrix {1 & \re\lam \\
0 & \im \lam }^\dag \bmatrix{\re (\Ls^H v x^H\Pb(\lam)x) \\
\im (\Ls^H v x^H\Pb(\lam)x)}$ when $\S= \herm,$ and  $\widehat{r} = r_s:= \bmatrix {1 & -\im\lam \\
0 & \re \lam }^\dag \bmatrix{\re (\Ls^H v x^H\Pb(\lam)x) \\
\im (\Ls^H v x^H\Pb(\lam)x)}$ when $\S= \skherm.$

Next we consider linearizations of $H$-even polynomials. Note that an $H$-even polynomial admits  $H$-even as well as $H$-odd linearizations and both have the same spectral symmetry as that of the polynomial. For purely imaginary eigenvalues of an $H$-even or $H$-odd polynomial, we have the following result.

\begin{theorem}\label{hheven2}
Let $\Pb \in \P$ be $H$-even and $(\lam, x) \in \C\times \C^n$ be such that $\|x\|_2 =1.$ Let $ \S \in \{\heven, \hodd\}$ and $ \Lb \in \S $ be a linearization of $\Pb$ with normalized ansatz vector $v.$ If $\lam \in i\R$ then  we have $$\sqrt{\frac{m+1}{2m}} \leq
\dfrac{\eta^{\S}_F(\lam,\Ls\otimes x,\Lb; v)}{\eta^{\heven}_F(\lam,x,\Pb)} \leq
\sqrt{2}  \mbox{ and } \sqrt{\frac{m+1}{2m}} \leq
\dfrac{\eta^{\S}_2(\lam,\Ls\otimes x,\Lb;
v)}{\eta^{\heven}_2(\lam,x,\Pb)} \leq 1.$$
The same bounds hold when $\Pb$ is $H$-odd.

\end{theorem}

\noin{\bf Proof:} The proof is exactly the same as that of Theorem~\ref{hherm2} and follows from Theorem~\ref{heven}. $\blacksquare$

This shows that a structured linearization of an $H$-even polynomial has least influence on the backward errors of approximate eigenelements when the approximate eigenvalues are purely imaginary. On the other hand, when the approximate eigenvalues are not purely imaginary, the structured backward errors of approximate eigenelements are not amenable to easy comparisons. Indeed, under the assumptions of Theorem~\ref{hheven2}, when $ \lam \in \C\setminus i\R$ by Theorem~\ref{heven}, we have
$$ \dfrac{\eta^\S_F(\lam,\Ls\otimes
x,\Lb; v)}{\eta(\lam, x, \Pb)} \leq \sqrt{2 + \frac{\|\widehat{r}\|_2^2}{\|r\|_2^2}} \mbox{ and }
\dfrac{\eta^\S_2(\lam,\Ls\otimes
x,\Lb; v)}{\eta(\lam, x, \Pb)} \leq \sqrt{1 + \frac{\|\widehat{r}\|_2^2}{\|r\|_2^2}},$$
where $ r := -\Pb(\lam)x$ and $ \widehat{r} = r_s$ when $ \S= \heven,$ and $  \widehat{r} = r_h$ when $ \S= \hodd.$

The obvious conclusion that we can draw is that real eigenvalues of a Hermitian/skew-Hermitian matrix polynomial can be computed either by solving a Hermitian  or a skew-Hermitian linearization. However, for non real eigenvalues it may be a good idea to solve Hermitian as well as skew-Hermitian linearizations and choose an eigenpair $(\lam, x)$ from Hermitian or skew-Hermitian linearization according as  $r_h \leq r_s$ or $ r_s \leq r_h.$ Similar conclusion holds for $H$-even/$H$-odd matrix polynomials. These observations are consistent with those made in~\cite{BaR09:cond} by analyzing the structured condition numbers of eigenvalues.

\section{Structured pseudospectra of structured polynomials}\label{chapter4:pseudo}
Let $\Pb $ be a regular polynomial. For $\lam \in \C,$ the
backward error of $\lam$ as an approximate eigenvalue of $\Pb$ is
given by $\eta(\lam, \Pb) := \min\{ \eta(\lam, x, \Pb) : x \in
\C^n \mbox{ and } \|x\|_2 =1\}.$ Since $\eta(\lam, x, \Pb) =
\|\Pb(\lam)x\|/{\|x\|_2\|\Lam\|_2},$ it follows that for the spectral as well as
for the Frobenius norms on $\C^{n\times n},$ we have
$$\eta(\lam, \Pb) := \frac{\sigmin(\Pb(\lam))}{\|\Lam\|_2}.$$
Similarly, for $ M \in \{2, F\}$ we define structured backward error of an approximate
eigenvalue $\lam$ of $\Pb \in \S$ by
$$\eta^{\S}_M(\lam, \Pb) := \min\{ \eta^{\S}_M(\lam, x, \Pb) : x \in
\C^n \mbox{ and } \|x\|_2 =1 \}.$$ In this section, we make an attempt to determine
$\eta^{\S}_M(\lam, \Pb).$ The backward errors of approximate
eigenvalues  and pseudospectra of a polynomial are closely
related. For $\ep
>0,$ the unstructured $\ep$-pseudospectrum of $\Pb$, denoted by $\elam(\Pb),$ is given by~(see~\cite{safique:thesis,sa2})
$$\sigma_\ep(\Pb) = \bigcup_{\nrm{\Delta\Pb}_M
\leq \ep}\{\sig(\Pb+\Delta\Pb): \Delta \Pb \in \P\} .$$  Obviously, we
have $\elam(\Pb) =\{ z \in \C : \eta(z, \Pb) \leq \ep \},$
assuming, for simplicity, that $\infty \notin \elam(\Pb),$ see~\cite{safique:thesis,sa2}.  For the
sake of simplicity in this section we make an implicit
assumption that $\infty \notin \elam(\Pb).$ Observe that since $\eta(\lam, \Pb)$ is the same for the spectral norm
and the Frobenius norm, we conclude that $\elam(\Pb)$ is the same for the spectral and the Frobenius norms.
Similarly, when $\Pb \in\S,$ we define the
structured $\ep$-pseudospectrum of $\Pb,$ denoted by
$\sig^{\S}_{\ep}(\Pb),$ by $$ \sig^{\S}_{\ep}(\Pb) :=
\bigcup_{\nrm{\Delta\Pb}_M \leq \ep}\{\sig(\Pb+\Delta\Pb) :
\Delta\Pb \in \S\}.$$ Then it follows that $\sig^{\S}_{\ep}(\Pb) =
\{ z \in \C : \eta^{\S}_M(\lam, \Pb) \leq \ep \}.$

\begin{theorem} Let $\S \in \{\sym, \sksym\}$
and $\Pb\in  \S.$ Then for the spectral norm, we have
$ \eta_2^{\S}(\lam, \Pb) = \eta(\lam, \Pb) \mbox{ and  }
\sig^{\S}_{\ep}(\Pb) = \sig_\ep(\Pb).$ On the other hand, for the Frobenius norm, we have
$ \eta_F^{\S}(\lam, \Pb) = \sqrt{2} \, \eta(\lam, \Pb)
\,\,\mbox{and}\,\, \sig^{\S}_{\ep}(\Pb) =
\sig_{\ep/\sqrt{2}}(\Pb)$ when $\S= \sksym,$  and
$ \eta_F^{\S}(\lam, \Pb) = \eta(\lam, \Pb) \mbox{ and  }
\sig^{\S}_{\ep}(\Pb) = \sig_\ep(\Pb)$ when $\S=\sym.$
\end{theorem}

\noin\pf For the spectral norm, by Theorem~\ref{tsym}, we have
$\eta^{\S}_2(\lam, x, \Pb) = \eta(\lam, x, \Pb)$ for all $x.$
Consequently, we have $ \eta^{\S}_2(\lam, \Pb) = \eta(\lam, \Pb).$
Hence the result follows.

For the Frobenius norm, the result follows from
Theorem~\ref{tskewsym} when $\Pb$ is skew-symmetric. So,
suppose that $\Pb$ is symmetric. Then $\Pb(\lam) \in
\C^{n\times n}$ is symmetric. Consider the Takagi factorization
$P(\lambda)=U\Sigma U^T,$ where $U$ is unitary and $\Sigma $ is a
diagonal matrix containing singular values of $\Pb(\lam)$ (appear
in descending order). Set $\sig := \Sigma(n,n)$ and $u := U(:,
n).$ Then we have $\Pb(\lam)\overline{u} = \sig u.$ Now define
$$\Delta A_j :=-\frac{\overline{\lam}^j \sig \, uu^T}{\|\Lam\|^2_2},$$
and consider the polynomial $\Delta\Pb(z) =\sum_{j=0}^m z^j
\Delta{A_j}.$ Then $\Delta\Pb$ is symmetric and
$\Pb(\lam)\overline{u}+\Delta\Pb(\lam) \overline{u} = 0.$ Note
that for $M \in \{2, F\}$  we have
$$\eta^{\S}_M(\lam, \Pb) \leq \nrm{\Delta\Pb}_M =
\frac{\sig}{\|\Lam\|_2} = \eta(\lam, \Pb) \mbox{ and hence }
\elam(\Pb) = \slam(\Pb).$$ This completes the proof.
$\blacksquare$

 When $\Pb$ is symmetric, the above proof shows how to
construct a symmetric polynomial $\Delta\Pb$ such that $\lam \in
\sig(\Pb+\Delta\Pb)$ and $\nrm{\Delta\Pb}_M = \eta^{\S}_M(\lam, \Pb).$
When $\Pb$ is skew-symmetric, using Takagi factorization of
the complex skew-symmetric matrix $\Pb(\lam),$ one can construct a
skew-symmetric polynomial $\Delta\Pb$ such that $\lam \in
\sig(\Pb+\Delta\Pb)$ and $\nrm{\Delta\Pb}_M = \eta^{\S}_M(\lam, \Pb).$
Indeed, consider the Takagi factorization $$ \Pb(\lam) =
U\diag(d_1, \cdots, d_m) U^T,$$ where $U$ is unitary, $d_j :=
\bmatrix{ 0 & s_j\\ -s_j & 0},$ $s_j \in \C$ is nonzero and
$|s_j|$ are singular values of $\Pb(\lam).$ Here the blocks $d_j$
appear in descending order of magnitude of $|s_j|.$ Note that
$\Pb(\lam) \overline{U} = U \diag(d_1, \cdots, d_m).$ Let $ u :=
U(:, n-1:n).$ Then $\Pb(\lam) \overline{u} = u d_m = u d_m
u^T\overline{u}.$ Now define
$$\Delta A_j := -\frac{\overline{\lambda^j} \, ud_m
u^T}{\|v_\lam\|_2^2}$$ and consider the pencil $\Delta\Pb(z) =
\sum_{j=0}^m z^j \Delta{A_j}.$ Then $\Delta\Pb$ is
skew-symmetric and $\Pb(\lam)\overline{u}+\Delta\Pb(\lam)
\overline{u} = 0.$ Further, we
have $$\eta^{\S}_2(\lam, \Pb) =\nrm{\Delta\Pb}_2 =
\frac{\sigmin(\Pb(\lam))}{\|\Lam\|_2} = \eta(\lam, \Pb),\,\,
\eta^{\S}_F(\lam, \Pb) = \nrm{\Delta\Pb}_F = \sqrt{2} \,
\frac{\sigmin(\Pb(\lam))}{\|\Lam\|_2} = \sqrt{2} \, \eta(\lam,
\Pb).$$

 We denote the unit circle in $\C$ by $\T,$ that $\T := \{ z
\in \C : |z| =1 \}.$ Then for the $T$-even or $T$-odd polynomials
we have the following result.

\begin{theorem} Let $\S \in \{T\mbox{-even, }  T\mbox{-odd}\}$ and $\Pb \in
\S$ and $m$ is odd.   Then for $\lam \in \T$ and the Frobenius
norm  we have
$ \eta^{\S}_F(\lam, \Pb) = \sqrt{2} \, \eta(\lam, \Pb) \mbox{ and  }
\sig^{\S}_{\ep}(\Pb)\cap \T = \sig_{\ep/{\sqrt{2}}}(\Pb)\cap \T.$
\end{theorem}

\noin\pf Let $ \lam \in \T.$ Then by Theorem~\ref{teven} and
Theorem~\ref{todd},  $ \eta^{\S}_F(\lam, x, \Pb) =
\sqrt{2}\, \|\Pb(\lam)x\|_2/{\|\Lam\|_2}$ for all $x$ such
that  $ \|x\|_2 =1.$ Hence taking minimum over $\|x\|_2 =1,$ we
obtain the desired results. $\blacksquare$

\begin{theorem}\label{pseudo} Let $\S \in \{ \herm, \skherm\}$ and
$\Pb \in \S.$  Then for $M \in \{2, F\}$ and $\lam \in \R,$  we have $ \eta^{\S}_M(\lam,
\Pb) = \eta(\lam, \Pb)$ and $ \sig^{\S}_{\ep}(\Pb)\cap \R =
\sig_{\ep}(\Pb)\cap \R.$

\end{theorem}

\noin\pf Note that $\Pb(\lam)$ is either Hermitian or
skew-Hermitian. Let $(\mu, u)$ be an eigenpair of the matrix
$\Pb(\lam)$ such that $|\mu| = \sigmin(\Pb(\lam))$ and $u^Hu=1.$
Then $\Pb(\lam)u = \mu u.$ Define $$\Delta{A_j} := -\frac{\lam^j
\, \mu \, uu^H}{\|\Lam\|_2^2}$$ and consider the polynomial
$\Delta\Pb(z) = \sum_{j=0}^m z^j \Delta{A_j}.$ Then $\Delta\Pb \in
\S$ and $\lam \in \sig(\Pb+\Delta\Pb).$ Further,  we have $\dm{\nrm{\Delta\Pb}_M =
\frac{\sigmin(\Pb(\lam))}{\|\Lam\|_2}.}$ Hence the result follows.
$\blacksquare$

\begin{theorem} Let $\S \in \{ H\mbox{-even, } H\mbox{-odd}\}$ and
$\Pb \in \S.$  Then for $M \in \{2, F\}$ and $\lam \in i\R,$  we have $ \eta^{\S}_M(\lam,
\Pb) = \eta(\lam, \Pb)$ and  $ \sig^{\S}_{\ep}(\Pb)\cap i
\,\R = \sig_{\ep}(\Pb)\cap i \,\R.$
\end{theorem}

\noin\pf Note for $\lam \in i \,\R,$ then the matrix $\Pb(\lam)$
is again either is Hermitian or skew-Hermitian. Hence the result
follows from the proof of Theorem~\ref{pseudo}. $\blacksquare$

We mention that the above results can be easily extended to
the case of general structured polynomials where the coefficients
matrices are elements of Jordan and/or Lie algebras.

Finally, we mention that the partial equality $\slam(\Lb)\cap
\Omega = \elam(\Lb)\cap \Omega,$ for an appropriate $\Omega \subset
\C,$ and the minimal perturbations constructed above as well as in section~3 are expected to be key tools for solution of certain structured distance problems, see~\cite{AlaBKMM}. For illustration, consider an $H$-even polynomial $\Pb.$ We have seen that the eigenvalues of $\Pb$  have Hamiltonian spectral symmetry, that is, the spectrum of $\Pb$ is symmetric with respect to the imaginary axis in the complex plane and thus appear in the pair $(\lam, - \overline{\lam}).$ Obviously the spectral symmetry degenerates if there are purely imaginary  eigenvalues. Often in practice the polynomial $\Pb$ is obtained by an approximation of the exact problem. Thus it may be the case that even though the exact problem has no purely imaginary eigenvalues but due to approximation error the polynomial $\Pb$ may have a few purely imaginary eigenvalues. So, the task is to construct a minimal perturbation $\Delta \Pb$ such that $\Pb+\Delta \Pb$  is $H$-even and that the perturbed polynomial $\Pb+\Delta \Pb$ has no eigenvalues on the imaginary axis. On the other hand, it may also be the case that all eigenvalues of the exact problem are purely imaginary (e.g. gyroscopic system) but due to approximation error the polynomial $\Pb$ has a few eigenvalues off the imaginary axis. In such a case the task is to construct a minimal  perturbation $\Delta\Pb$ such that $\Pb+\Delta \Pb$ is $H$-even and that all eigenvalues of $\Pb+\Delta \Pb$ are purely imaginary. The minimal perturbations so constructed and the partial equality of between structured and unstructured pseudospectra so established are expected to be key tools in solving these problems and will be investigated elsewhere.

\von
\noin{\bf Conclusion.} We have derived backward errors of approximate eigenelements of structured matrix polynomials. We have constructed minimal structured perturbations that induce the approximate eigenelements as exact eigenelements of the perturbed polynomials. The minimal perturbations so constructed are expected to be useful in analyzing the evolution of eigenvalues of structured polynomials under structure preserving perturbations. We have analyzed the influence of structure preserving linearizations on the approximate eigenelements of structured matrix polynomials. Also, we have provided a recipe for selecting structure preserving linearizations so that the linearizations have almost no adverse effect on the approximate eigenelements of the structured matrix polynomials. We have briefly analyzed structured pseudospectra of structured matrix polynomials and have shown that partial equality between structured and unstructured pseudospectra holds for certain structured polynomials. These results are expected to be useful for constructing minimal structured perturbations that move eigenvalues of the structured polynomials along certain directions which in turn are expected to be key tools for solving certain structured distance problems.

\end{document}